\documentclass[11pt,a4paper,draft]{article}

\usepackage[reqno]{amsmath}
\usepackage{amssymb,amsthm,upref,amscd}
\usepackage[T1]{fontenc}
\usepackage{times}

\newtheorem{theorem}{Theorem}[section]
\newtheorem{corollary}[theorem]{Corollary}

\newtheorem{lemma}[theorem]{Lemma}
\newtheorem{proposition}[theorem]{Proposition}
\newtheorem{remark}[theorem]{Remark}

\theoremstyle{definition} \theoremstyle{remark}
\numberwithin{equation}{section}

\begin{document}

\title{\textbf{Two Positive Normalized Solutions and Phase Separation for Coupled Schr\"{o}dinger Equations on Bounded Domain with $L^2$-Supercritical and Sobolev Critical or Subcritical Exponent\ \\
}}
\author{Linjie Song$^{\mathrm{a,}}$\thanks{%
Song is supported by "Shuimu Tsinghua Scholar Program". Email: songlinjie@mail.tsinghua.edu.cn.}\
\;\; and \;\ Wenming Zou$^{\mathrm{a}}$\thanks{%
Zou is supported by NSFC(12171265). Email: zou-wm@mail.tsinghua.edu.cn.} \ \\
\\
{\small $^{\mathrm{a}}$Department of Mathematical Sciences, Tsinghua University, Beijing 100084, China}
}
\date{}
\maketitle

\begin{abstract}
In this paper we study the existence of positive normalized solutions of the following coupled Schr\"{o}dinger system:
\begin{align}
\left\{
\begin{aligned}
& -\Delta u = \lambda_u u + \mu_1 u^3 + \beta uv^2, \quad x \in \Omega, \\
& -\Delta v = \lambda_v v + \mu_2 v^3 + \beta u^2 v, \quad x \in \Omega, \\
& u > 0, v > 0 \quad \text{in } \Omega, \quad u = v = 0 \quad \text{on } \partial\Omega,
\end{aligned}
\right. \nonumber
\end{align}
with the $L^2$ constraint
\begin{align}
\int_{\Omega}|u|^2dx = c_1, \quad \quad \int_{\Omega}|v|^2dx = c_2, \nonumber
\end{align}
where $\mu_1, \mu_2 > 0$, $\beta \neq 0$, $c_1, c_2 > 0$, and $\Omega \subset \mathbb{R}^N$ ($N = 3, 4$) is smooth, bounded, and star-shaped. Note that the nonlinearities and the coupling terms are both $L^2$-supercritical in dimensions 3 and 4, Sobolev subcritical in dimension 3, Sobolev critical in dimension 4. We show that this system has a positive normalized solution which is a local minimizer. We further show that the system has a second positive normalized solution, which is of M-P type when $N = 3$. This seems to be the first existence result of two positive normalized solutions for such a Schr\"{o}dinger system, especially in the Sobolev critical case. We also study the limit behavior of the positive normalized solutions in the repulsive case $\beta \to -\infty$, and phase separation is expected.
\bigskip



\noindent\textbf{Data availability statement:} The manuscript has no associate data.

\end{abstract}

\medskip

\section{Introduction}

Our main goal in this paper is to establish the existence of two positive solutions for the following elliptic system:
\begin{equation} \label{eq1.1}
\left\{
\begin{aligned}
& -\Delta u = \lambda_u u + \mu_1 u^3 + \beta uv^2, \quad x \in \Omega, \\
& -\Delta v = \lambda_v v + \mu_2 v^3 + \beta u^2 v, \quad x \in \Omega, \\
& u > 0, v > 0 \quad \text{in } \Omega, \quad u = v = 0 \quad \text{on } \partial\Omega,
\end{aligned}
\right.
\end{equation}
with the $L^2$ constraint
\begin{align} \label{eq1.2}
\int_{\Omega}|u|^2dx = c_1, \quad \quad \int_{\Omega}|v|^2dx = c_2,
\end{align}
where $\mu_1, \mu_2 > 0$ and $\beta \neq 0$ is a coupling constant, $c_1, c_2 > 0$, and $\Omega \subset \mathbb{R}^N$ ($N = 3, 4$) is smooth, bounded, and star-shaped with respect to the origin.

\vskip0.1in
System \eqref{eq1.1} comes from finding standing wave solutions of time-dependent coupled
nonlinear Schr\"{o}dinger systems:
\begin{align} \label{eq1.5}
\left\{
\begin{aligned}
& i\partial_{t} \Phi_1 + \Delta \Phi_1 + \mu_1|\Phi_1|^2\Phi_1 + \beta |\Phi_2|^2\Phi_1= 0, \quad x \in \Omega, \ t > 0, \\
& i\partial_{t} \Phi_2 + \Delta \Phi_2 + \mu_2|\Phi_2|^2\Phi_2 + \beta |\Phi_1|^2\Phi_2= 0, \quad x \in \Omega, \ t > 0, \\
& \Phi_j = \Phi_j(t,x) \in \mathbb{C}, \quad  j =1,2 \\
& \Phi_j(t,x) = 0, \quad x \in \partial \Omega, \quad t > 0, \ j =1,2.
\end{aligned}
\right.
\end{align}
When $N \leq 3$, system \eqref{eq1.1} appears in many physical problems, especially in nonlinear optics. Physically, the solution $\Phi_j$ denotes the $j$th component of the beam in Kerr-like photorefractive media (see \cite{Akhmediev-Ankiewicz}). The positive constant $\mu_j$ is for self-focusing in the $j$th component of the beam. The coupling constant $\beta$ is the interaction between the two components of the beam. The problem \eqref{eq1.1} also arises in the Hartree-Fock theory for a double condensate, that is, a binary mixture of Bose-Einstein condensates in two different hyperfine states $|1\rangle$ and $|2\rangle$ (see \cite{EGBB}). Further, $\Phi_j$ values are the corresponding condensate amplitudes, $\mu_j$ and $\beta$ are the intraspecies and interspecies scattering lengths. The sign of $\beta$ determines whether the interactions of states $|1\rangle$ and $|2\rangle$ are repulsive or attractive: the interaction is attractive if $\beta > 0$, and the interaction is repulsive if $\beta < 0$, where the two states are in strong competition.

\vskip0.1in

In previous literatures, people often found solutions of \eqref{eq1.1} with $\lambda_u, \lambda_v$ fixed and without the $L^2$ constraints, see \cite{Ambrosetti-Colorado-2006,Ambrosetti-Colorado-2007,Bartsch-Wang,Bartsch-Wang-Wei,Lin-Wei,Maia-Pellacci-Squassina,Sirakov} for $N \leq 3$ and see \cite{Chen-Zou-2012} for $N = 4$. What we aim to establish in this paper are \textit{normalized solutions}, i.e., solutions with the constraints that $\int_{\Omega}|u|^2dx = c_1, \int_{\Omega}|v|^2dx = c_2$. This is more suitable for physical settings. Moreover, establishing the existence of normalized solutions is crucial to study orbital stability/instability for standing waves. Also, people can easily see that any normalized solution is fully nontrivial, i.e., $u \neq 0, v \neq 0$. Solutions on entire space have been widely studied and the $L^2$-critical exponent $2 + 4/N$ plays an important role, see e.g. \cite{BLZ,BJN,BN,LZ} and their references. For a single equation, in a series work \cite{HPS,HS,HS2,Song,Song2,Song3,Song4}, the first author with his co-authors developed new frameworks to prove the \textit{existence} and \textit{uniqueness} of normalized solutions. Furthermore, people can refer \cite{CHMS,HLS} for biharmonic NLS and for \cite{HLS,EHLS} on the waveguide manifold. Particularly, in \cite{Song5} which is in preparation, the authors aim to establish the existence of two normalized positive solutions for a semi-linear elliptic equation with general nonlinearities on star-shaped bounded domains. As for normalized solutions for systems on bounded domains, the only paper we know is \cite{NTV2}, in which one positive solution was found. To the best of our knowledge, there is no existence result of two normalized solutions for \eqref{eq1.1}. We also remark that obtaining a second solution is much more difficult than the first one, especially when $N = 4$ (the Sobolev critical case). Our discussions can be extended to higher dimensional situations ($N \geq 5$).

\vskip0.1in

Let $\lambda_1 := \lambda_1(\Omega)$ be the first eigenvalue of $-\Delta$ on $\Omega$ with Dirichlet boundary condition and $e_1$ be the corresponding positive unit eigenfunction.  The first step is to find a \textit{local minimizer}. The structure of a local minimum has been observed in \cite{NTV2}. We will provide the result and its proof for readers' convenience though the proof is not totally new. To find solutions of (\ref{eq1.1}), we search for critical points of the energy
\begin{equation}
E(u,v) = \frac{1}{2}\int_{\Omega}(|\nabla u|^{2}+|\nabla v|^{2})dx - \frac{1}{4}\int_{\Omega}(\mu_1|u^+|^4+2\beta|u^+|^2|v^+|^2+\mu_2|v^+|^4)dx,
\end{equation}
under the constraint
\begin{equation}
\int_{\Omega}|u^+|^{2}dx = c_1, \int_{\Omega}|v^+|^{2}dx = c_2,
\end{equation}
where $u^+ = \max\{u,0\}$, $v^+ = \max\{v,0\}$. Let $\vec{c} = (c_1,c_2)$. We further set
$$
\mathcal{H} := H_0^1(\Omega) \times H_0^1(\Omega),
$$
$$
S_{\vec{c}}^+ := \left\{(u,v) \in \mathcal{H}: \int_{\Omega}|u^+|^2dx = c_1, \int_{\Omega}|v^+|^2dx = c_2\right\},
$$
$$
\mathcal{B}_\rho := \left\{(u,v) \in \mathcal{H}: \int_{\Omega}(|\nabla u|^{2}+|\nabla v|^{2})dx < \rho^2\right\}.
$$
Let $\mathcal{C}_N$ be the Sobolev best constant of $H_0^1(\Omega) \hookrightarrow L^4(\Omega)$ where $\Omega \subset \mathbb{R}^N$, i.e.,
\begin{equation}
 	\mathcal{C}_N =\sup_{u \in H_0^1(\Omega) \setminus \{0\} } \cfrac{\left(\int_{\Omega} |u|^{4}dx\right)^{\frac{1}{4}}}{\left(\int_{\Omega} |\nabla u|^2 dx\right)^{\frac{1}{2}}}.
\end{equation}
Note that $\mathcal{C}_4 = 1/\sqrt{\mathcal{S}}$, where
$\mathcal{S}$ is the Sobolev best constant of $\mathcal{D}^{1,2}(\mathbb{R}^4)\hookrightarrow L^{4}(\mathbb{R}^4)$,
\begin{equation}\label{sobolev constant}
 	\mathcal{S}=\inf_{u \in \mathcal{D}^{1,2}(\mathbb{R}^4) \setminus \{0\} } \cfrac{\int_{\mathbb{R}^4} |\nabla u|^2 dx}{\left(\int_{\mathbb{R}^4} |u|^{4}dx\right)^{\frac{1}{2}}},
\end{equation}
where $\mathcal{D}^{1,2}(\mathbb{R}^4)=\{u\in L^{4}(\mathbb{R}^4): |\nabla u| \in L^{2}(\mathbb{R}^4)\}$ with norm $\left\| u \right\|_{\mathcal{D}^{1,2}}:=\left(\int_{\mathbb{R}^4}|\nabla u|^2dx\right)^{\frac{1}{2}} $.

\vskip0.1in
Now we state our theorems on the existence of positive normalized solutions.

\begin{theorem}
\label{thmB.2}
Let $\Omega$ be Lipschitz and bounded. When $\beta < 0$, let $c_1, c_2$ satisfy
\begin{align} \label{eq < 0}
\lambda_1(c_1+c_2-\beta\mathcal{C}_N c_1c_2) \leq \frac{1}{2\mathcal{C}_N\max\{\mu_1,\mu_2\}}.
\end{align}
When $\beta > 0$, let $c_1, c_2$ satisfy
\begin{align}
\lambda_1(c_1+c_2) \leq \frac{1}{\mathcal{C}_N(\sqrt{(\mu_1-\mu_2)+4\beta^2}+\mu_1+\mu_2)}.
\end{align}
Then the value
$$
\nu_{\vec{c},\rho}:= \inf_{S_{\vec{c}}^+ \cap \overline{\mathcal{B}_\rho}}E,
$$
is achieved in $\mathcal{B}_\rho$, for a suitable $\rho > \sqrt{\lambda_1(c_1+c_2)}$. Furthermore, there exists a critical point $(u_{\vec{c}},v_{\vec{c}}) \in S_{\vec{c}}^+$ such that, for some $\lambda_{u_{\vec{c}}}$, $\lambda_{v_{\vec{c}}}$,
\begin{equation} \label{eq of u}
\left\{
\begin{aligned}
&-\Delta u_{\vec{c}} = \lambda_{u_{\vec{c}}} u_{\vec{c}} + \mu_1 u_{\vec{c}}^3 + \beta u_{\vec{c}}v_{\vec{c}}^2, \quad x \in \Omega, \\
&-\Delta v_{\vec{c}} = \lambda_{v_{\vec{c}}} v_{\vec{c}} + \mu_2 v_{\vec{c}}^3 + \beta u_{\vec{c}}^2 v_{\vec{c}}, \quad x \in \Omega, \\
\end{aligned}
\right.
\end{equation}
with $\displaystyle \int_\Omega(|\nabla u_{\vec{c}}|^2+|\nabla v_{\vec{c}}|^2)dx < \rho^2$, and $u_{\vec{c}} > 0$, $v_{\vec{c}} > 0$.
\end{theorem}

\begin{remark}
In above theorem to obtain a local minimizer, it is not necessary to assume that $\Omega$ is star-shaped. Such an assumption will be used to get the second solution, particularly in the step using the Pohozaev identity to construct a bounded (PS) sequence.
\end{remark}

Next we aim to obtain the second positive normalized solution. To do this, we will construct a M-P structure on $S_{\vec{c}}^+$.  The main difficulty in Sobolev subcritical case ($N = 3$) is to find a \textit{bounded (PS) sequence} and we further face the challenge of losing the \textit{compactness} in Sobolev critical case ($N = 4$). Pohozaev manifold method is often used to get a bounded (PS) sequence in entire space cases. This depends heavily on the scaling invariance of $\mathbb{R}^N$, and is not applicable to star-shaped bounded domain cases. To overcome this difficulty, our  idea is to use the Pohozaev identity. We will use the monotonicity trick to find a (PS) sequence and then use the Pohozaev identity to show its boundedness. To overcome the difficulty of losing compactness when $N = 4$, we further estimate the value of the M-P level. Compared with classical Br\'ezis-Nirenberg question, we construct a path under the $L^2$ constraint, that, combined with the existence of a local minimizer at an energy level below the mountain pass level, makes things more complex. Similar difficulties were encountered in \cite{JL,WW} for a Sobolev critical Schr\"{o}dinger equation on the entire space. Since $\Omega$ is not scaling invariant, difficulties become harder to overcome and we need to develop some new techniques.

\begin{theorem}
\label{thmB.5} Let $\Omega$ be smooth, bounded, and star-shaped with respect to the origin. Under the assumptions of Theorem \ref{thmB.2}, the system \eqref{eq1.1}-\eqref{eq1.2} has the second positive solution $(\tilde{u}_{\vec{c}},\tilde{v}_{\vec{c}}) \neq (u_{\vec{c}},v_{\vec{c}})$. Furthermore, if $N = 3$, then $(\tilde{u}_{\vec{c}},\tilde{v}_{\vec{c}})$ is of M-P type.
\end{theorem}

\begin{remark}
	In \cite{Chen-Zou-2012}, the second author with Z.J. Chen studied the existence of positive least energy solution for
	\begin{equation} \label{eqcz}
	\left\{
	\begin{aligned}
	& -\Delta u = \lambda_u u + \mu_1 u^3 + \beta uv^2, \quad x \in \Omega \subset \mathbb{R}^4, \\
	& -\Delta v = \lambda_v v + \mu_2 v^3 + \beta u^2 v, \quad x \in \Omega \subset \mathbb{R}^4, \\
	& u > 0, v > 0 \quad \text{in } \Omega, \quad u = v = 0 \quad \text{on } \partial \Omega, \\
	\end{aligned}
	\right.
	\end{equation}
	where $\mu_1, \mu_2 > 0$, $\lambda_u, \lambda_v \in (0,\lambda_1(\Omega))$ are fixed, without the $L^2$ constraint. They showed an existence result when $\beta \in (-\infty,0)\cup(0,\beta_1)\cup(\beta_2,\infty)$, where $\beta_1 \in (0,\min\{\mu_1,\mu_2\}]$ and $\beta_2 \geq \max\{\mu_1,\mu_2\}$. We remark that the condition on $\beta$ is not needed when searching for positive normalized solutions with the $L^2$ constraint. Indeed, as the above Theorems \ref{thmB.2} and \ref{thmB.5} have shown, our results hold true for all $\beta \neq 0$.
\end{remark}

\begin{remark}
	In \cite{Chen-Zou-2012}, a nonexistence result was provided. In fact, multiply the first equation of \eqref{eq1.1} by $v$, the second equation of \eqref{eq1.1} by $u$, and integrate over $\Omega$, which yields
	$$
	(\lambda_v - \lambda_u)\int_\Omega uvdx + (\beta - \mu_1)\int_\Omega u^3vdx + (\mu_2 - \beta)\int_\Omega uv^3dx = 0.
	$$
	This shows \eqref{eqcz} does not have a nontrivial nonnegative solution when $\lambda_u < \lambda_v$ and $\mu_1 \leq \beta \leq \mu_2$. It is open whether \eqref{eqcz} has a (fully) nontrivial positive solution for some $\lambda_u \geq \lambda_v$ when $\mu_1 \leq \beta \leq \mu_2$. Our Theorems \ref{thmB.2} and \ref{thmB.5} give a positive answer for this open question.
\end{remark}

Next, we study the limit behavior of the positive normalized solutions in the repulsive case $\beta \to -\infty$. Note that $\beta$ tending to $-\infty$ is not allowed in \eqref{eq < 0} for fixed $c_1 > 0, c_2 > 0$. Fortunately, by Remark \ref{rmk2.2} below we can provide a condition on $c_1, c_2$ independent of $\beta < 0$ such that two positive normalized solutions exist. It is expected that components of the limiting profile tend to separate in different regions of the underlying domain $\Omega$. This phenomenon, called \textit{phase separation}, has been well studied for cases when $\lambda_u, \lambda_v$ are fixed in \cite{Chen-Zou-2012,NTTV,WW1,WW2} and for other kinds of elliptic systems with strong competition in \cite{CL,CR,CTV}. Let $\lambda_2 := \lambda_2(\Omega)$ be the second eigenvalue of $-\Delta$ on $\Omega$ with Dirichlet boundary condition and $e_2$ be the corresponding unit eigenfunction. Then we have the following result.

\begin{theorem} \label{thmps}
	Let $\Omega$ be smooth, bounded, and star-shaped with respect to the origin and $\beta < 0$. Let $c_1 > 0, c_2 > 0$ satisfy
	\begin{align} \label{eq beta < 0}
	\lambda_2(c_1+c_2) \leq \frac{1}{2\mathcal{C}_N\max\{\mu_1,\mu_2\}}.
	\end{align}
	Then \eqref{eq1.1}-\eqref{eq1.2} has two positive solutions $(u_{1,\beta},v_{1,\beta})$ (a local minimizer) and $(u_{2,\beta},v_{2,\beta})$, with Lagrange multipliers $(\eta_{1,\beta},\xi_{1,\beta})$ and $(\eta_{2,\beta},\xi_{2,\beta})$ respectively. $\{(u_{1,\beta},v_{1,\beta})\}_{\beta < 0}$, $\{(u_{2,\beta},v_{2,\beta})\}_{\beta < 0}$ are uniformly bounded in $H_0^1(\Omega) \times H_0^1(\Omega)$ and $C^{0,\alpha}(\overline{\Omega}) \times C^{0,\alpha}(\overline{\Omega})$ for $\alpha \in (0,1)$. Further, up to subsequences, $u_{i,\beta} \to Q_i^+$, $v_{i,\beta} \to -Q_i^-$ in $H_0^1(\Omega) \cap C^{0,\alpha}(\overline{\Omega})$ as $\beta \to -\infty$, where $Q_i$ solves
	\begin{equation} \label{eqlimit}
	\left\{
	\begin{aligned}
	&-\Delta Q_i = \eta_i Q_i^+ + \xi_i Q_i^- + \mu_1(Q_i^+)^3 + \mu_2(Q_i^-)^3, \quad x \in \Omega, \\
    &\int_{\Omega}|Q_i^+|^2dx = c_1, \quad \int_{\Omega}|Q_i^-|^2dx = c_2, \quad Q_i \in H_0^1(\Omega),
	\end{aligned}
	\right.
	\end{equation}
	where $Q_i^+ = \max\{Q_i,0\}$, $Q_i^- = \min\{Q_i,0\}$, $\eta_i = \displaystyle\lim_{\beta \to -\infty}\eta_{i,\beta}$, and $\xi_i = \displaystyle\lim_{\beta \to -\infty}\xi_{i,\beta}$, $i = 1,2$.
\end{theorem}

\begin{remark}
	The results in Theorem \ref{thmps} for $(u_{1,\beta},v_{1,\beta})$ were proved in \cite{NTV2}. For $(u_{2,\beta},v_{2,\beta})$, borrowing an idea in \cite[Section 5]{NTV2}, we can complete the proof. We remark that the proof for $(u_{2,\beta},v_{2,\beta})$ is not trivial since we need new ideas to give uniform $\displaystyle  H^1_0(\Omega)$ bounds for $\{u_{2,\beta}\}_{\beta < 0}$ and $\displaystyle \{v_{2,\beta}\}_{\beta < 0}$. For $(u_{1,\beta},v_{1,\beta})$ this is trivial since $(u_{1,\beta},v_{1,\beta}) \in \mathcal{B}_\rho$ where $\rho$ is independent of $\beta < 0$. For $\displaystyle  (u_{2,\beta},v_{2,\beta})$, we first give a uniform upper bound for $\displaystyle \{E(u_{2,\beta},v_{2,\beta})\}_{\beta < 0}$ (see \eqref{eq bound} below) and then use the Pohozaev identity to provide uniform $H^1_0(\Omega)$ bounds, see more details in Section \ref{secps}.
\end{remark}

The paper is organized as follows. In Section \ref{local minimizer} we provide the proof for Theorem \ref{thmB.2}. In Section \ref{boundedness} we establish the existence of a bounded (PS) sequence with some properties. In Section \ref{estimate} we obtain the estimation of the M-P level in the Sobolev critical case. In Section \ref{second solution} we complete the proof of Theorem \ref{thmB.5}. Finally in Section \ref{secps} we prove Theorem \ref{thmps}.

\

\noindent\textbf{Acknowledgement. }{\small \textit{ The first author thanks T.H. Liu for his careful reading and suggestions to address the case when $\beta < 0$.} }

\section{Proof for Theorem \ref{thmB.2}} \label{local minimizer}

We firstly state a result which will be used in the proof of Theorem \ref{thmB.2}. We will prove it at the end of this section.

\begin{lemma}
\label{thmB.4} Under the assumptions of Theorem \ref{thmB.2}, we have
\begin{align} \label{eq pre for m-p of E}
E(\sqrt{c_1}e_1,\sqrt{c_2}e_1) < \inf_{(u,v) \in \partial\mathcal{B}_{\rho}\cap S_{\vec{c}}^+}E(u,v)
\end{align}
for some $\rho > \sqrt{\lambda_1(c_1+c_2)}$. If $N = 4$, we can further take $\rho$ such that
\begin{align}
& \rho^2 \leq \frac{2\mathcal{S}^2}{\max\{\mu_1,\mu_2\}}  \ \text{if} \ \beta < 0, \label{eq r1} \\
& \rho^2 \leq \frac{4\mathcal{S}^2}{\sqrt{(\mu_1-\mu_2)+4\beta^2}+\mu_1+\mu_2} \ \text{if} \ \beta > 0. \label{eq r2}
\end{align}
\end{lemma}

\begin{remark} \label{rmk2.2}
As readers will see, \eqref{eq pre for m-p of E} is used to show that
\begin{align} \label{eq ext}
\inf_{(u,v) \in \mathcal{B}_{\rho}\cap S_{\vec{c}}^+}E(u,v) < \inf_{(u,v) \in \partial\mathcal{B}_{\rho}\cap S_{\vec{c}}^+}E(u,v).
\end{align}
We note that \eqref{eq pre for m-p of E} can be replaced by other ones, which can give other ranges of $c_1, c_2$ such that the positive normalized solutions exist. Indeed, let $\lambda_2 := \lambda_2(\Omega)$ be the second eigenvalue of $-\Delta$ on $\Omega$ with Dirichlet boundary condition and $e_2$ be the corresponding unit eigenfunction. Note that $e_2$ is sign-changing. Let $e_2^+ = \max\{e_2,0\}$ and $e_2^- = e_2 - e_2^+$. We assume that
$$
E\left(\sqrt{c_1}\frac{e_2^+}{\|e_2^+\|_{L^2(\Omega)}}, -\sqrt{c_2}\frac{e_2^-}{\|e_2^-\|_{L^2(\Omega)}}\right) < \inf_{(u,v) \in \partial\mathcal{B}_{\rho}\cap S_{\vec{c}}^+}E(u,v).
$$
Then \eqref{eq ext} holds true for $\rho > \sqrt{\lambda_2(c_1+c_2)}$.
\end{remark}

\subsection{Proof of Theorem \ref{thmB.2} assuming Lemma \ref{thmB.4}}

Before proving Theorem \ref{thmB.2}, we show the following result.

\begin{theorem}
\label{thmB.3} Let $\rho^2 > \lambda_1(c_1+c_2)$, $\vec{c} = (c_1,c_2)$, $c_1 > 0$, $c_2 > 0$ and $\tau > 0$ be fixed, and $\nu_{\vec{c},\rho}$ be defined by Theorem \ref{thmB.2}. If $N = 4$, we further assume that $\rho^{2^{\ast}-2} \leq \frac{2^{\ast}}{2}S^{2^{\ast}}$. If
\begin{equation} \label{eqB.2}
\nu_{\vec{c},\rho} < \widehat{\nu}_{\vec{c},\rho} := \inf_{S_{\vec{c}}^+ \cap \overline{\mathcal{B}_\rho} \setminus \mathcal{B}_{\rho - \tau}}E,
\end{equation}
then $E|_{S_{\vec{c}}^+}$ contains a critical point at the level $\nu_{\vec{c},\rho}$ contained in $\mathcal{B}_{\rho}$.
\end{theorem}

\noindent{\bf Proof.} By H\"{o}lder inequality and Sobolev embedding, one gets $-\infty < \nu_{\vec{c},\rho}, \widehat{\nu}_{\vec{c},\rho} < +\infty$. By the assumption that $\nu_{\vec{c},\rho} < \widehat{\nu}_{\vec{c},\rho}$, we can take a minimizing sequence $\{(u_n,v_n)\} \subset \mathcal{B}_{\rho - \tau}$, implying that $\{(u_n,v_n)\}$ is bounded in $\mathcal{H}$. Using Ekeland variational principle, we may assume that $E|_{S_{\vec{c}}^+}'(u_n,v_n) \to 0$ as $n \to \infty$. Up to a subsequence, assume that $u_n \rightharpoonup u, v_n \rightharpoonup v$ weakly in $H_0^1(\Omega)$ and $u_n \to u, v_n \to v$ strongly in $L^2(\Omega)$.

If $N = 3$, we further assume that $u_n \to u, v_n \to v$ strongly in $L^4(\Omega)$. The strong convergences of $u_n, v_n$ in $L^2(\Omega)$ yield that $(u,v) \in S_{\vec{c}}^+$. Further, by the weakly lower semi continuity of the norm, one knows that
\begin{align} \label{eq rho}
\int_{\Omega}|\nabla u|^2dx \leq \liminf_{n \to \infty}\int_{\Omega}|\nabla u_n|^2dx, \quad \int_{\Omega}|\nabla v|^2dx \leq \liminf_{n \to \infty}\int_{\Omega}|\nabla v_n|^2dx.
\end{align}
Hence $(u,v) \in \mathcal{B}_{\rho - \tau}$. On the one hand, using $(u,v) \in S_{\vec{c}}^+ \cap \mathcal{B}_{\rho - \tau}$ we have $E(u,v) \geq \nu_{\vec{c},\rho}$. On the other hand, by \eqref{eq rho} and strong convergences of $u_n, v_n$ in $L^4(\Omega)$, one gets $E(u,v) \leq \displaystyle \lim_{n \to \infty}E(u_n,v_n) = \nu_{\vec{c},\rho}$. Thus $E(u,v) = \nu_{\vec{c},\rho}$. This in turn shows that
\begin{align}
\int_{\Omega}|\nabla u|^2dx = \liminf_{n \to \infty}\int_{\Omega}|\nabla u_n|^2dx, \quad \int_{\Omega}|\nabla v|^2dx = \liminf_{n \to \infty}\int_{\Omega}|\nabla v_n|^2dx.
\end{align}
Note that $u_n \rightharpoonup u, v_n \rightharpoonup v$ weakly in $H_0^1(\Omega)$. Thus $u_n \to u, v_n \to v$ strongly in $H_0^1(\Omega)$. This completes the proof when $N = 3$.

If $N = 4$, we further assume that
\begin{align}
\begin{aligned}
u_n \rightharpoonup u, \ \quad v_n \rightharpoonup v, \ \ \quad & \text{weakly in} \ L^4(\Omega), \nonumber \\
u_n^2 \rightharpoonup u^2, \quad v_n^2 \rightharpoonup v^2, \quad & \text{weakly in} \ L^2(\Omega), \nonumber \\
u_n^3 \rightharpoonup u^3, \quad  v_n^3 \rightharpoonup v^3, \quad & \text{weakly in} \ L^{4/3}(\Omega), \nonumber \\
u_n \to u, \ \quad v_n \to v, \ \ \quad & \text{strongly in} \ L^2(\Omega), \nonumber \\
u_n \to u, \ \quad v_n \to v, \ \ \quad & \text{for almost everywhere} \ x \in \Omega. \nonumber
\end{aligned}
\end{align}
Let $z_n = u_n - u$, $w_n = v_n - v$. Then by Br\'ezis-Lieb Lemma (see \cite{BL} and \cite[formula 5.37]{Chen-Zou-2012}), we have
\begin{align}
& \int_{\Omega}|u_n^+|^4dx = \int_{\Omega}|u^+|^4dx + \int_{\Omega}|z_n^+|^4dx + o_n(1), \\
& \int_{\Omega}|v_n^+|^4dx = \int_{\Omega}|v^+|^4dx + \int_{\Omega}|w_n^+|^4dx + o_n(1), \\
& \int_{\Omega}|u_n^+|^2|v_n^+|^2dx = \int_{\Omega}|u^+|^2|v^+|^2dx + \int_{\Omega}|z_n^+|^2|w_n^+|^2dx + o_n(1).
\end{align}
Hence, we deduce that
\begin{align} \label{eq2.11}
E(u_n,v_n) = E(u,v) + E(z_n,w_n) + o_n(1).
\end{align}
Similar to the case of $N = 3$, $(u,v) \in S_{\vec{c}}^+ \cap \mathcal{B}_\rho$ and thus $E(u,v) \geq \nu_{\vec{c},\rho}$. Then using \eqref{eq2.11} we obtain $E(z_n,w_n) \leq o_n(1)$. Up to a subsequence, assume that $\int_{\Omega}|\nabla z_n|^2dx \to k \geq 0, \int_{\Omega}|\nabla w_n|^2dx \to l \geq 0$ as $n \to \infty$. When $\beta < 0$, $E(z_n,w_n) \leq o_n(1)$ implies that
\begin{align}
(k + l) & \leq \lim_{n \to \infty}\frac{1}{2}\int_{\Omega}(\mu_1|z_n^+|^{4}+\mu_2|w_n^+|^{4})dx \nonumber \\
& \leq \frac{1}{2\mathcal{S}^2}(\mu_1k^2 + \mu_2l^2) \nonumber \\
& \leq \frac{\max\{\mu_1,\mu_2\}}{2\mathcal{S}^2}(k+l)^2.
\end{align}
This yields either $k + l = 0$ or $k+l \geq \frac{2\mathcal{S}^2}{\max\{\mu_1,\mu_2\}}$.
From
\begin{align}
& \int_{\Omega}|\nabla u_n|^2dx = \int_{\Omega}|\nabla u|^2dx + \int_{\Omega}|\nabla z_n|^2dx + o_n(1), \nonumber \\
& \int_{\Omega}|\nabla v_n|^2dx = \int_{\Omega}|\nabla v|^2dx + \int_{\Omega}|\nabla w_n|^2dx + o_n(1), \nonumber
\end{align}
we derive that $\int_{\Omega}(|\nabla z_n|^2+|\nabla w_n|^2)dx < (\rho-\tau)^2$, and so $k + l < \rho^2$. Since $\rho^2 \leq \frac{2\mathcal{S}^2}{\max\{\mu_1,\mu_2\}}$, we obtain $k + l = 0$ and this implies $k = 0, l = 0$. One gets $u_n \to u, v_n \to v$ strongly in $H_0^1(\Omega)$.

When $\beta > 0$, we have $2\beta a b \leq \beta\epsilon a^2 + \beta/\epsilon b^2$ for any $a > 0, b > 0$. Take
$$
\epsilon = \frac{\sqrt{(\mu_1-\mu_2)^2+4\beta^2}-\mu_1+\mu_2}{2\beta}.
$$
Then there holds
\begin{align} \label{eq esti when > 0}
& \int_{\Omega}(\mu_1|z_n^+|^4+2\beta|z_n^+|^2|w_n^+|^2+\mu_2|w_n^+|^4)dx \nonumber \\
\leq & \frac{\sqrt{(\mu_1-\mu_2)+4\beta^2}+\mu_1+\mu_2}{2}\int_{\Omega}(|z_n^+|^4+|w_n^+|^4)dx.
\end{align}
Using $E(z_n,w_n) \leq o_n(1)$ we obtain that
\begin{align}
(k + l) & \leq \lim_{n \to \infty}\frac{1}{2}\int_{\Omega}(\mu_1|z_n^+|^4+2\beta|z_n^+|^2|w_n^+|^2+\mu_2|w_n^+|^4)dx \nonumber \\
& \leq \frac{\sqrt{(\mu_1-\mu_2)+4\beta^2}+\mu_1+\mu_2}{4\mathcal{S}^2}(k^2 + l^2) \nonumber \\
& \leq \frac{\sqrt{(\mu_1-\mu_2)+4\beta^2}+\mu_1+\mu_2}{4\mathcal{S}^2}(k+l)^2.
\end{align}
This yields either $k + l = 0$ or $k+l \geq \frac{4\mathcal{S}^2}{\sqrt{(\mu_1-\mu_2)+4\beta^2}+\mu_1+\mu_2}$.
Similar to the case $\beta < 0$, we have $k + l < \rho^2$. Since $\rho^2 \leq \frac{4\mathcal{S}^2}{\sqrt{(\mu_1-\mu_2)+4\beta^2}+\mu_1+\mu_2}$, we obtain $k + l = 0$ and this implies $k = 0, l = 0$. One gets $u_n \to u, v_n \to v$ strongly in $H_0^1(\Omega)$. The proof is complete now.
\qed\vskip 10pt

\noindent \textbf{Proof of Theorem \ref{thmB.2} assuming Lemma \ref{thmB.4}.  } Since $\rho^2 > \lambda_1(c_1+c_2)$, we obtain that $(\sqrt{c_1}e_1,\sqrt{c_2}e_1) \in \mathcal{B}_\rho$, implying that $\nu_{\vec{c},\rho} \leq E(\sqrt{c_1}e_1,\sqrt{c_2}e_1)$. Using Lemma \ref{thmB.4}, it holds
$$
E(\sqrt{c_1}e_1,\sqrt{c_2}e_1) < \inf_{S_{\vec{c}}^+ \cap \partial\mathcal{B}_\rho}E = \inf_{S_{\vec{c}}^+ \cap \overline{\mathcal{B}_\rho} \setminus \mathcal{B}_{\rho - \tau}}E + o_\tau(1),
$$
where $o_\tau(1) \to 0$ as $\tau \to 0$. Hence, we can derive that \eqref{eqB.2} holds for small $\tau > 0$.

Then by Theorem \ref{thmB.3}, $E|_{S_{\vec{c}}^+}$ contains a critical point $(u_{\vec{c}},v_{\vec{c}})$ at the level $\nu_{\vec{c},\rho}$ contained in $\mathcal{B}_{\rho}$. By Lagrange multiplier principle, $(u_{\vec{c}},v_{\vec{c}})$ satisfies
\begin{equation}
\left\{
\begin{aligned}
&-\Delta u_{\vec{c}} = \lambda_{u_{\vec{c}}} u_{\vec{c}}^+ + \mu_1 (u_{\vec{c}}^+)^3 + \beta u_{\vec{c}}^+(v_{\vec{c}}^+)^2, \quad x \in \Omega, \\
&-\Delta v_{\vec{c}} = \lambda_{v_{\vec{c}}} v_{\vec{c}}^+ + \mu_2 (v_{\vec{c}}^+)^3 + \beta (u_{\vec{c}}^+)^2 v_{\vec{c}}^+, \quad x \in \Omega,
\end{aligned}
\right. \nonumber
\end{equation}
for some $\lambda_{u_{\vec{c}}}$, $\lambda_{v_{\vec{c}}}$. Multiplying $u_{\vec{c}}^-$, $v_{\vec{c}}^-$ and integrating on $\Omega$ for these two equations respectively, we obtain $$\int_\Omega|\nabla u_{\vec{c}}^-|^2dx = \int_\Omega|\nabla v_{\vec{c}}^-|^2dx = 0,$$ implying that $u_{\vec{c}}^- = v_{\vec{c}}^- = 0$ and hence $u_{\vec{c}}, v_{\vec{c}} \geq 0$. By strong maximum principle, $u_{\vec{c}}, v_{\vec{c}} > 0$. Thus, $\int_\Omega|u_{\vec{c}}|^2dx = \int_\Omega|u_{\vec{c}}^+|^2dx = c_1$, $\int_\Omega|v_{\vec{c}}|^2dx = \int_\Omega|v_{\vec{c}}^+|^2dx = c_2$ (hence $(u_{\vec{c}},v_{\vec{c}})$ is fully nontrivial) and $(u_{\vec{c}},v_{\vec{c}})$ solves \eqref{eq of u}. The proof is complete.
\qed\vskip 5pt

\subsection{Proof of Lemma \ref{thmB.4}}

\textbf{Proof of Lemma \ref{thmB.4}.  } When $\beta < 0$, we have
\begin{align}
\inf_{(u,v) \in \partial\mathcal{B}_\rho\cap S_{\vec{c}}^+}E(u,v) & \geq \frac12\rho^2 - \frac{\mathcal{C}_N\max\{\mu_1,\mu_2\}}{4}\rho^4 \nonumber \\
& = -\frac{\mathcal{C}_N\max\{\mu_1,\mu_2\}}{4}\left(\rho^2 - \frac{1}{\mathcal{C}_N\max\{\mu_1,\mu_2\}}\right) \nonumber \\
&\quad+ \frac{1}{4\mathcal{C}_N\max\{\mu_1,\mu_2\}}.
\end{align}
Take $$\rho^2 = \frac{1}{\mathcal{C}_N\max\{\mu_1,\mu_2\}} > \lambda_1(c_1+c_2).$$ Also, note that $\rho^2 < \frac{2\mathcal{S}^2}{\max\{\mu_1,\mu_2\}}$. We have that
\begin{align} \label{eq < 0 2}
\inf_{(u,v) \in \partial\mathcal{B}_\rho\cap S_{\vec{c}}^+}E(u,v) \geq \frac{1}{4\mathcal{C}_N\max\{\mu_1,\mu_2\}}.
\end{align}
Since
$$
-2\beta c_1c_2\int_{\Omega}e_1^4dx \leq -2\beta c_1c_2\mathcal{C}_N\lambda_1,
$$
we have
$$
E(\sqrt{c_1}e_1,\sqrt{c_2}e_1) < \frac{\lambda_1(c_1+c_2-\beta\mathcal{C}_N c_1c_2)}{2}.
$$
Using \eqref{eq < 0} and \eqref{eq < 0 2}, one gets
$$
E(\sqrt{c_1}e_1,\sqrt{c_2}e_1) < \inf_{(u,v) \in \partial\mathcal{B}_\rho\cap S_{\vec{c}}^+}E(u,v).
$$

When $\beta > 0$, using \eqref{eq esti when > 0} we can get
\begin{align}
& \inf_{(u,v) \in \partial\mathcal{B}_\rho\cap S_{\vec{c}}^+}E(u,v) \nonumber \\
\geq & \frac12\rho^2 -\frac{\mathcal{C}_N(\sqrt{(\mu_1-\mu_2)+4\beta^2}+\mu_1+\mu_2)}{8}\rho^4 \nonumber \\
= & -\frac{\mathcal{C}_N(\sqrt{(\mu_1-\mu_2)+4\beta^2}+\mu_1+\mu_2)}{8}\left( \rho^2-\frac{2}{\mathcal{C}_N(\sqrt{(\mu_1-\mu_2)+4\beta^2}+\mu_1+\mu_2}\right) \nonumber \\
& + \frac{1}{2\mathcal{C}_N(\sqrt{(\mu_1-\mu_2)+4\beta^2}+\mu_1+\mu_2)}.
\end{align}
Take $$\rho^2 = \frac{2}{\mathcal{C}_N(\sqrt{(\mu_1-\mu_2)+4\beta^2}+\mu_1+\mu_2)} > \lambda_1(c_1+c_2).$$ Also, note that $\rho^2 < \frac{4\mathcal{S}^2}{\sqrt{(\mu_1-\mu_2)+4\beta^2}+\mu_1+\mu_2}$. We have that
\begin{align}
\inf_{(u,v) \in \partial\mathcal{B}_\rho\cap S_{\vec{c}}^+}E(u,v) & \geq \frac{1}{2\mathcal{C}_N(\sqrt{(\mu_1-\mu_2)+4\beta^2}+\mu_1+\mu_2)} \nonumber \\
& \geq \frac{\lambda_1(c_1+c_2)}{2} \nonumber \\
& > E(\sqrt{c_1}e_1,\sqrt{c_2}e_1).
\end{align}

\qed\vskip 5pt

\section{The bounded (PS) sequence} \label{boundedness}

From now on, we focus on the existence of the second positive solution. To obtain a bounded (PS) sequence, we introduce the family of functionals
$$E_\theta(u,v) =$$
\begin{align}
\left\{
	\begin{aligned}
	\frac{1}{2}\int_{\Omega}(|\nabla u|^{2}+|\nabla v|^{2})dx - \frac{\theta}{4}\int_{\Omega}(\mu_1|u^+|^4+2\beta|u^+|^2|v^+|^2+\mu_2|v^+|^4)dx, \beta > 0, \\
	\frac{1}{2}\int_{\Omega}(|\nabla u|^{2}+|\nabla v|^{2}-\beta|u^+|^2|v^+|^2)dx - \frac{\theta}{4}\int_{\Omega}(\mu_1|u^+|^4+\mu_2|v^+|^4)dx,  \beta < 0,
	\end{aligned}
	\right. \nonumber
\end{align}
where $\theta \in [1/2,1]$. The crucial idea is to use the monotonicity trick \cite{BCJS,Jean} and the Pohozaev identity.

\begin{lemma}[Uniform M-P geometry] \label{M-P geometry}
Assume that \eqref{eq pre for m-p of E} holds for some $\rho > \sqrt{(c_1+c_2)\lambda_1}$, then there exists $\epsilon \in (0,1/2)$ such that
\begin{align} \label{eq pre for m-p}
E_\theta(\sqrt{c_1}e_1,\sqrt{c_2}e_1) + \delta < \inf_{(u,v) \in \partial\mathcal{B}_{\rho}\cap S_{\vec{c}}^+}E_\theta(u,v)
\end{align}
holds for some $\delta > 0$ independent of $\theta \in [1-\epsilon,1]$.
Furthermore, there exists $(z,w) \in S_{\vec{c}}^+ \cap \mathcal{B}_{\rho}^c$ such that
\begin{align}
m_\theta := \inf_{\gamma \in \Gamma}\sup_{t \in [0,1]}E_\theta(\gamma(t)) > & E_\theta(\sqrt{c_1}e_1,\sqrt{c_2}e_1) + \delta \nonumber \\
= & \max\{E_\theta(\sqrt{c_1}e_1,\sqrt{c_2}e_1),E_\theta(z,w)\} + \delta
\end{align}
where
$$
\Gamma := \{\gamma \in C([0,1],S_{\vec{c}}^+): \gamma(0) = (\sqrt{c_1}e_1,\sqrt{c_2}e_1), \gamma(1) = (z,w)\}
$$
is independent of $\theta$.
\end{lemma}

\noindent\textbf{Proof.  }
Take $2\delta = \displaystyle\inf_{(u,v) \in \partial\mathcal{B}_{\rho}\cap S_{\vec{c}}^+}E(u,v) - E(\sqrt{c_1}e_1,\sqrt{c_2}e_1)$. Then
\begin{align}
& E_\theta(\sqrt{c_1}e_1,\sqrt{c_2}e_1) \nonumber \\
& = E(\sqrt{c}e_1,\sqrt{c_2}e_1) + \frac{1-\theta}{4}(\mu_1c_1^2 + 2\beta c_1c_2 +\mu_2c_2^2) \int_\Omega e_1^4dx \nonumber \\
& = \inf_{(u,v) \in \partial\mathcal{B}_{\rho}\cap S_{\vec{c}}^+}E(u,v) - 2\delta + O(1-\theta) \nonumber \\
& = \inf_{(u,v) \in \partial\mathcal{B}_{\rho}\cap S_{\vec{c}}^+}\left( E_\theta(u,v) + \frac{\theta-1}{4}\int_\Omega(\mu_1|u^+|^4+2\beta|u^+|^2|v^+|^2+\mu_2|v^+|^4)dx\right) \nonumber \\
& ~~~~~~ - 2\delta + O(1-\theta) \nonumber \\
& = \inf_{(u,v) \in \partial\mathcal{B}_{\rho}\cap S_{\vec{c}}^+}E_\theta(u,v) - 2\delta + O(1-\theta).
\end{align}
By choosing $\epsilon \in (0,1/2)$ such that $|O(1-\theta)| < \delta$ for all $\theta \in [1-\epsilon,1]$, we have
$$
E_\theta(\sqrt{c_1}e_1,\sqrt{c_2}e_1) + \delta < \inf_{(u,v) \in \partial\mathcal{B}_{\rho}\cap S_{\vec{c}}^+}E_\theta(u,v), \quad \forall \theta \in [1-\epsilon,1].
$$
Next, for $(\phi,\psi) \in S_{\vec{c}}^+$ where $\phi,\;\;  \psi$ have compact supports and $\text{supp}\phi \cap \text{supp}\psi = \emptyset$, we set
$$
\phi_t = t^{\frac{N}{2}}\phi(tx),\;\;  \psi_t = t^{\frac{N}{2}}\psi(tx), t \geq 1.
$$
Then $(\phi_t,\psi_t) \in S_{\vec{c}}^+$ and
\begin{align}
&E_\theta(\phi_t, \psi_t)\nonumber \\
& = \frac{1}{2}\int_{\Omega}(|\nabla \phi_t|^{2} + |\nabla \psi_t|^{2})dx - \frac{\theta}{4}\int_{\Omega}(\mu_1|\phi_t^+|^4+2\beta|\phi_t^+|^2|\psi_t^+|^2+\mu_2|\psi_t^+|^4)dx \nonumber \\
&\leq \frac{t^2}{2}\int_{\Omega}(|\nabla \phi|^{2} + |\nabla \psi|^{2})dx - \frac{t^{N}}{8}\int_{\Omega}(\mu_1|\phi^+|^4+\mu_2|\psi^+|^4)dx \nonumber \\
&\to  -\infty,
\end{align}
as $t \to \infty$ since $N > 2$.

Take $(z,w) = (\phi_t,\psi_t)$ with $t$ large enough such that $E_\theta(z,w) < E_\theta(\sqrt{c_1}e_1,\sqrt{c_2}e_1)$. Also, we have $\|\nabla z\|_{L^2(\Omega)}^2 + \|\nabla w\|_{L^2(\Omega)}^2 > \rho^2$. Then for any $\gamma \in \Gamma$, there exists $t^\ast \in (0,1)$ such that $\gamma(t^\ast) \in \partial\mathcal{B}_{\rho}$. Hence, $$\inf_{\gamma \in \Gamma}\sup_{t \in [0,1]}E_\theta(\gamma(t)) \geq \inf_{(u,v) \in \partial\mathcal{B}_{\rho} \cap S_{\vec{c}}^+}E_\theta(u,v).$$ By \eqref{eq pre for m-p} we complete the proof.
\qed\vskip 5pt

\begin{lemma} \label{lem equal of m}
Under the assumptions of Lemma \ref{M-P geometry}, let $\mathcal{P}_\theta \subset S_{\vec{c}}^+$ be the maximal path connected branch containing $(\sqrt{c_1}e_1,\sqrt{c_2}e_1)$ such that
$$
E_\theta(u,v) < m_\theta, \forall (u,v) \in \mathcal{P}_\theta.
$$
Define
\begin{align}
\tilde{m}_\theta := \inf_{\gamma \in \tilde{\Gamma}}\sup_{t \in [0,1]}E_\theta(\gamma(t)) > E_\theta(u_0,v_0) = \max\{E_\theta(u_0,v_0),E_\theta(z,w)\}
\end{align}
where $(u_0,v_0) \in \mathcal{P}_\theta$ is a positive solution of \eqref{eq1.1}, $(z,w)$ is given by Lemma \ref{M-P geometry},
$$
\tilde{\Gamma} := \left\{\gamma \in C([0,1],S_{\vec{c}}^+): \gamma(0) = (u_0,v_0), \gamma(1) = (z,w)\right\}.
$$
Then $\tilde{m}_\theta = m_\theta$.
\end{lemma}

\begin{remark}
	Before proving Lemma \ref{lem equal of m}, we give an equivalent characteristic of $\mathcal{P}_\theta$ for a better understanding of this definition. Let
	$$
	\Sigma_\theta = \left\{\sigma \in C([0,1],S_{\vec{c}}^+): \gamma(0) = (\sqrt{c_1}e_1,\sqrt{c_2}e_1), \sup_{t \in [0,1]}E_\theta(\sigma(t)) < m_\theta\right\}.
	$$
	On the one hand, $\sigma(1) \in \mathcal{P}_\theta$ for all $\sigma \in \Sigma_\theta$. On the other hand, for any $(u,v) \in \mathcal{P}_\theta$, there exists a path $\sigma \in \Sigma_\theta$ such that $(u,v) = \sigma(1)$. Hence, $\mathcal{P}_\theta$ is the image of the following mapping:
	$$
	\Sigma_\theta \to S_{\vec{c}}^+, \quad \sigma \mapsto \sigma(1),
	$$
	i.e.,
	$$
	\mathcal{P}_\theta = \{\sigma(1): \sigma \in \Sigma_\theta\}.
	$$
	We introduce $\mathcal{P}_\theta$ for estimating $m_\theta$. As readers will see, we can provide an estimation for $\tilde{m}_\theta$ firstly (see Theorem \ref{Estimate of the M-P level} below), then using $\tilde{m}_\theta = m_\theta$ we obtain the estimation of $m_\theta$.
\end{remark}

\noindent\textbf{Proof.}  By the definition of $\mathcal{P}_\theta$ and the fact that $(u_0,v_0) \in \mathcal{P}_\theta$, we can take $\gamma_0 \in C([0,1],\mathcal{P}_\theta)$ such that $\gamma_0(0) = (\sqrt{c_1}e_1, \sqrt{c_2}e_1), \gamma_0(1) = (u_0,v_0)$. For any $\tilde{\gamma} \in \tilde{\Gamma}$, we define
\begin{equation}
\gamma(t) =
\left\{
\begin{aligned}
&\gamma_0(2t), & 0 \leq t < \frac12, \\
&\tilde{\gamma}(2t -1), & \frac12 < t \leq 1.
\end{aligned}
\right.
\end{equation}
Then $\gamma \in \Gamma$ and hence $\displaystyle\sup_{t \in [0,1]}E_\theta(\gamma(t)) \geq m_\theta$. From the definition of $\mathcal{P}_\theta$, $\displaystyle\sup_{t \in [0,1]}E_\theta(\gamma_0(t)) < m_\theta$. Therefore, $\displaystyle\sup_{t \in [0,1]}E_\theta(\tilde{\gamma}(t)) \geq m_\theta$, implying that $\tilde{m}_\theta \geq m_\theta$.

On the contrary, for any $\gamma \in \Gamma$, we define
\begin{equation}
\tilde{\gamma}(t) =
\left\{
\begin{aligned}
& \gamma_0(1-2t), & 0 \leq t < \frac12, \\
& \gamma(2t -1), & \frac12 < t \leq 1.
\end{aligned}
\right.
\end{equation}
Then $\tilde{\gamma} \in \tilde{\Gamma}$ and hence $\displaystyle\sup_{t \in [0,1]}E_\theta(\tilde{\gamma}(t)) \geq \tilde{m}_\theta$. From the definition of $\mathcal{P}_\theta$, $\displaystyle\sup_{t \in [0,1]}E_\theta(\gamma_0(t)) < m_\theta \leq \tilde{m}_\theta$. Therefore, $\displaystyle\sup_{t \in [0,1]}E_\theta(\gamma(t)) \geq \tilde{m}_\theta$, implying that $\tilde{m}_\theta \leq m_\theta$. The proof is complete.
\qed\vskip 5pt

At this point we wish to use the monotonicity trick \cite{BCJS,Jean} on the family of functionals $E_\theta$. We recall the general setting in which the theorem is stated. Let $(W,\langle\cdot,\cdot\rangle,\|\cdot\|)$ and $(H,(\cdot,\cdot),|\cdot|)$ be two Hilbert spaces, which form a variational triple, i.e.,
$$
W \hookrightarrow H = H' \hookrightarrow W'
$$
with continuous injections. For simplicity, we assume that the continuous injection $W \hookrightarrow H$ has the norm at most 1 and identify $W$ with its image in $H$. Let $P \subset H$ be a cone. Any $u \in H$ can be decomposed into $u^+ + u^-$ where $u^+ \in P$ and $u^- \in -P$. Define
$$
S_{\vec{c}}^+ := \{(u,v) \in W \times W: |u^+|^2 =c_1, |v^+|^2 =c_2\}, c_1, c_2 > 0.
$$
In our application there hold that $W = H_0^1(\Omega)$, $H = L^2(\Omega)$ and $P = \{u \in L^2(\Omega): u \geq 0\}$.
Using the results in \cite{BCJS,Jean}, we have the following theorem.

\begin{theorem}[Monotonicity trick] \label{monotonicity trick}
Let $I \subset \mathbb{R}^+$ be an interval. We consider a family $(E_\theta)_{\theta \in I}$ of $C^1$-functionals on $W \times W$ of the form
$$
E_\theta(u,v) = A(u,v) - \theta B(u,v), \theta \in I,
$$
where $B(u,v) \geq 0, \forall (u,v) \in W \times W$ and such that either $A(u,v) \rightarrow \infty$ or $B(u,v) \rightarrow \infty$ as $\|u\| + \|v\| \rightarrow \infty$. We assume there are two points $(z_1, w_1)$, $(z_2,w_2)$ in $S_{\vec{c}}^+$ (independent of $\theta$) such that setting
$$
\Gamma = \{\gamma \in C([0, 1], S_{\vec{c}}^+), \gamma(0) = (z_1, w_1), \gamma(1) = (z_2,w_2)\},
$$
there holds, $\forall \tau \in I$,
$$
m_\theta:= \inf_{\gamma \in \Gamma}\sup_{t \in [0,1]}E_\theta(\gamma(t)) > \max\{E_\theta(z_1, w_1),E_\theta(z_2,w_2)\}.
$$
Then, for almost every $\theta \in I$, there is a sequence $\{(\phi_n,\psi_n)\} \subset S_{\vec{c}}^+$ such that

$(i)$ $\{(\phi_n,\psi_n)\}$ is bounded in $W \times W$, $(ii)$ $E_\theta(\phi_n,\psi_n) \rightarrow m_\theta$, $(iii)$ $E'_\theta|_{S_{\vec{c}}^+}(\phi_n,\psi_n) \rightarrow 0$ in $W'\times W'$.
\end{theorem}

\begin{proposition}[Positive solutions for almost everywhere $\theta$]
\label{Positive solutions for almost every theta}
    Under the assumptions of Lemma \ref{M-P geometry}, for almost everywhere $\theta \in [1-\epsilon,1]$ where $\epsilon$ is given by Lemma \ref{M-P geometry}, there exists a critical point $(u_\theta,v_\theta)$ of $E_\theta$ constrained on $S_{\vec{c}}^+$, which solves
\begin{equation} \label{eq theta}
\left\{
\begin{aligned}
& -\Delta u_\theta = \lambda_{u_\theta} u_\theta + \theta \mu_1 u_\theta^3 + \theta\beta u_\theta v_\theta^2, \quad x \in \Omega, \\
& -\Delta v_\theta = \lambda_{v_\theta} v_\theta + \theta\mu_2 v_\theta^3 + \theta\beta u_\theta^2 v_\theta, \quad x \in \Omega, \\
& u_\theta > 0, v_\theta > 0 \quad \text{in } \Omega, \quad u_\theta = v_\theta = 0 \quad \text{on } \partial \Omega,
\end{aligned}
\right.
\end{equation}
when $\beta > 0$ and
\begin{equation} \label{eq theta 2}
\left\{
\begin{aligned}
& -\Delta u_\theta = \lambda_{u_\theta} u_\theta + \theta \mu_1 u_\theta^3 + \beta u_\theta v_\theta^2, \quad x \in \Omega, \\
& -\Delta v_\theta = \lambda_{v_\theta} v_\theta + \theta\mu_2 v_\theta^3 + \beta u_\theta^2 v_\theta, \quad x \in \Omega, \\
& u_\theta > 0, v_\theta > 0 \quad \text{in } \Omega, \quad u_\theta = v_\theta = 0 \quad \text{on } \partial \Omega,
\end{aligned}
\right.
\end{equation}
when $\beta < 0$, for some $\lambda_{u_\theta}, \lambda_{v_\theta}$. If $N = 3$, then $u_\theta$ is at the level $m_\theta$. If $N = 4$, then one of the following holds true:
\begin{itemize}
  \item[(a)] let $\mathcal{P}_\theta$ be given by Lemma \ref{lem equal of m}, then $(u_\theta,v_\theta) \in \partial\mathcal{P}_\theta$ and thus $E_\theta(u_\theta,v_\theta) = m_\theta$;
  \item[(b)] let $\mathcal{P}_\theta$ be given by Lemma \ref{lem equal of m}, then $(u_\theta,v_\theta) \notin \overline{\mathcal{P}_\theta}$ and $E_\theta(u_\theta,v_\theta) \leq m_\theta$.
\end{itemize}
Furthermore, if $\beta > 0$,
\begin{align} \label{inequality}
\int_{\Omega}(|\nabla u_\theta|^2+|\nabla v_\theta|^2)dx > \frac{\theta N}{4}\int_{\Omega}(\mu_1 |u_\theta|^4 + 2\beta |u_\theta|^2 |v_\theta|^2 + \mu_2 |v_\theta|^4)dx;
\end{align}
if $\beta < 0$,
\begin{align} \label{inequality 2}
\int_{\Omega}(|\nabla u_\theta|^2+|\nabla v_\theta|^2)dx > \frac{\beta N}2\int_{\Omega} |u_\theta|^2 |v_\theta|^2dx + \frac{\theta N}{4}\int_{\Omega}(\mu_1 |u_\theta|^4 + \mu_2 |v_\theta|^4)dx.
\end{align}
\end{proposition}

\noindent\textbf{Proof.} Applying Theorem \ref{monotonicity trick} when $\beta > 0$ with
$$
W = H_0^1(\Omega), \quad H = L^2(\Omega), \quad P = \{u \in L^2(\Omega): u \geq 0\},
$$
$$
A(u,v) = \frac{1}{2}\int_{\Omega}(|\nabla u|^{2}+|\nabla v|^{2})dx,
$$
$$
B(u,v) = \frac{1}{4}\int_{\Omega}(\mu_1|u^+|^4+2\beta|u^+|^2|v^+|^2+\mu_2|v^+|^4)dx;
$$
and when $\beta < 0$ with $$
W = H_0^1(\Omega), \quad H = L^2(\Omega), \quad P = \{u \in L^2(\Omega): u \geq 0\},
$$
$$
A(u,v) = \frac{1}{2}\int_{\Omega}(|\nabla u|^{2}+|\nabla v|^{2}-\beta|u^+|^2|v^+|^2)dx,
$$
$$
B(u,v) = \frac{1}{4}\int_{\Omega}(\mu_1|u^+|^4+\mu_2|v^+|^4)dx.
$$
By Lemma \ref{M-P geometry}, for almost everywhere $\theta \in [1-\epsilon,1]$, there exists a bounded (PS) sequence $\{(u_n,v_n)\} \subset S_{\vec{c}}^+$ satisfying $E_\theta(u_n,v_n) \rightarrow m_\theta$ and $E'_\theta|_{S_{\vec{c}}^+}(u_n,v_n) \rightarrow 0$ as $n \to \infty$.

If $N = 3$, using Sobolev compact embedding, through a standard process, one gets that $u_n \to u_\theta$, $v_n \to v_\theta$ strongly in $H_0^1(\Omega)$ passing to a subsequence if necessary. Hence, $(u_\theta,v_\theta)$ is a critical point of $E_\theta$ constrained on $S_{\vec{c}}^+$ at the level of $m_\theta$.

If $N = 4$, we need additional discussions. Up to a subsequence, we assume that
\begin{equation}
\begin{aligned}
&u_n \rightharpoonup u_\theta, \quad v_n \rightharpoonup v_\theta, \ \text{ weakly in } H_0^1(\Omega), \\
&u_n \rightharpoonup u_\theta, \quad v_n \rightharpoonup v_\theta, \ \text{ weakly in } L^{4}(\Omega),\\
&u_n \to u_\theta, \quad v_n \to v_\theta, \ \text{ strongly in } L^p(\Omega) \text{ for } 1 \leq p <4,\\
&u_n \to u_\theta, \quad v_n \to v_\theta, \ \text{ almost everywhere in } \Omega.
\end{aligned}
\end{equation}
It can be verified that $(u_\theta,v_\theta) \in S_{\vec{c}}^+$ is a critical point of $E_\theta$ constrained on $S_{\vec{c}}^+$. We claim that $(u_\theta,v_\theta) \notin \mathcal{P}_\theta$ where $\mathcal{P}_\theta$ is given by Lemma \ref{lem equal of m}. Suppose on the contrary that $(u_\theta,v_\theta) \in \mathcal{P}_\theta$. Then similar discussions to Theorem \ref{Estimate of the M-P level} in the next section yield that
\begin{align}
& m_\theta < E_\theta(u_\theta,v_\theta) + \frac{1}{4\theta\max\{\mu_1,\mu_2\}}\mathcal{S}^2 \ \text{if} \ \beta < 0, \nonumber \\
& m_\theta < E_\theta(u_\theta,v_\theta) + \frac{1}{2\theta(\sqrt{(\mu_1-\mu_2)^2+4\beta^2}+\mu_1+\mu_2)}\mathcal{S}^2 \ \text{if} \ \beta > 0.
\end{align}
Let $z_n = u_n - u_\theta$, $w_n = v_n - v_\theta$. Since $E'_\theta|_{S_{\vec{c}}^+}(u_n,v_n) \rightarrow 0$ as $n \to \infty$, there exists $\lambda_{n,1}$, $\lambda_{n,2}$ such that $\partial_uE_\theta(u_n,v_n) - \lambda_{n,1}u_n^+ \to 0$, $\partial_vE_\theta(u_n,v_n) - \lambda_{n,2}v_n^+ \to 0$, as $n \to \infty$. Let $\lambda_{u_\theta}$, $\lambda_{v_\theta}$ be the Lagrange multipliers correspond to $u_\theta$, $v_\theta$ respectively. Then for some $\psi \in H_0^1(\Omega)$ with $\int_\Omega u_\theta^+\psi dx \neq 0$, we obtain
$$
\lambda_{n,1} = \frac{1}{\int_\Omega u_n^+\psi dx}(\langle \partial_uE_\theta(u_n,v_n),\psi\rangle + o_n(1))$$
 $$\xrightarrow{n \to \infty} \frac{1}{\int_\Omega u_\theta^+\psi dx}\langle \partial_uE_\theta(u_\theta,v_\theta),\psi\rangle = \lambda_{u_\theta}.
$$
Similarly, we have $\lambda_{n,2} \to \lambda_{v_\theta}$ as $n \to \infty$. Using the Br\'{e}zis-Lieb Lemma (see \cite{BL}) and the facts that
$$
\partial_uE_\theta(u_n,v_n) - \lambda_{n,1}u_n^+ \to 0, \ \partial_uE_\theta(u_\theta,v_\theta) - \lambda_{u_\theta} u_\theta = 0,
$$
$$
\partial_vE_\theta(u_n,v_n) - \lambda_{n,2}v_n^+ \to 0, \ \partial_vE_\theta(u_\theta,v_\theta) - \lambda_{v_\theta} v_\theta = 0,
$$
when $\beta > 0$, one gets
\begin{align}
& \int_{\Omega}|\nabla z_n|^2dx = \theta\int_{\Omega}(\mu_1|z_n|^{4}+\beta|z_n|^{2}|w_n|^{2})dx + o_n(1), \label{eq of zn} \\
& \int_{\Omega}|\nabla w_n|^2dx = \theta\int_{\Omega}(\mu_2|w_n|^{4}+\beta|z_n|^{2}|w_n|^{2})dx + o_n(1). \label{eq of wn}
\end{align}
From the definition of $\mathcal{S}$ and \eqref{eq of zn}, \eqref{eq of wn} we deduce that
\begin{align*}
& \int_{\Omega}(|\nabla z_n|^2+|\nabla w_n|^2)dx + o_n(1) \\
= & \theta\int_{\Omega}(\mu_1|z_n|^{4}+2\beta|z_n|^{2}|w_n|^{2}+\mu_2|w_n|^{4})dx \\
\leq & \frac{\theta(\sqrt{(\mu_1-\mu_2)^2+4\beta^2}+\mu_1+\mu_2)}{2}\int_{\Omega}(|z_n|^{4}+|w_n|^{4})dx \\
\leq & \frac{\theta(\sqrt{(\mu_1-\mu_2)^2+4\beta^2}+\mu_1+\mu_2)}{2\mathcal{S}^2}\left( \left( \int_{\Omega}|\nabla z_n|^2dx\right) ^2 + \left( \int_{\Omega}|\nabla w_n|^2dx\right) ^2\right)  \\
\leq & \frac{\theta(\sqrt{(\mu_1-\mu_2)^2+4\beta^2}+\mu_1+\mu_2)}{2\mathcal{S}^2}\left(\int_{\Omega}(|\nabla z_n|^2+|\nabla w_n|^2)dx\right)^2.
\end{align*}
Let $\int_{\Omega}|\nabla z_n|^2dx \to k \geq 0, \int_{\Omega}|\nabla w_n|^2dx \to l \geq 0$ as $n \to \infty$. We have
\begin{align}
k + l \leq \frac{\theta(\sqrt{(\mu_1-\mu_2)^2+4\beta^2}+\mu_1+\mu_2)}{2\mathcal{S}^2}(k+l)^2.
\end{align}
If $k + l > 0$, one gets $k+l \geq \frac{2\mathcal{S}^2}{\theta(\sqrt{(\mu_1-\mu_2)^2+4\beta^2}+\mu_1+\mu_2)}$.
This shows that
\begin{align}
E_\theta(z_n,w_n) & = \frac{1}{4}(k + l) + o_n(1) \nonumber \\
& \geq \frac{1}{2\theta(\sqrt{(\mu_1-\mu_2)^2+4\beta^2}+\mu_1+\mu_2)}\mathcal{S}^2 + o_n(1).
\end{align}
Then, using the Br\'{e}zis-Lieb Lemma (see \cite{BL}) one gets
\begin{align}
E_\theta(u_n,v_n) & = E_\theta(u_\theta,v_\theta) + E_\theta(z_n,w_n) + o_n(1) \nonumber \\
& \geq E_\theta(u_\theta,v_\theta) + \frac{1}{2\theta(\sqrt{(\mu_1-\mu_2)^2+4\beta^2}+\mu_1+\mu_2)}\mathcal{S}^2 + o_n(1),
\end{align}
which contradicts to
$$
E_\theta(u_n,v_n) \to m_\theta < E_\theta(u_\theta,v_\theta) + \frac{1}{2\theta(\sqrt{(\mu_1-\mu_2)^2+4\beta^2}+\mu_1+\mu_2)}\mathcal{S}^2.
$$
Thus $k = l = 0$.

Similar to \eqref{eq of zn} and \eqref{eq of wn}, when $\beta < 0$, one gets
\begin{align}
& \int_{\Omega}|\nabla z_n|^2dx = \theta\int_{\Omega}\mu_1|z_n|^{4}dx + \int_{\Omega}\beta|z_n|^{2}|w_n|^{2}dx + o_n(1), \label{eq of zn 2} \\
& \int_{\Omega}|\nabla w_n|^2dx = \theta\int_{\Omega}\mu_2|w_n|^{4}dx + \int_{\Omega}\beta|z_n|^{2}|w_n|^{2}dx + o_n(1). \label{eq of wn 2}
\end{align}
From the definition of $\mathcal{S}$ and \eqref{eq of zn 2}, \eqref{eq of wn 2} we deduce that
\begin{align*}
& ~~~~~~ \int_{\Omega}(|\nabla z_n|^2+|\nabla w_n|^2)dx + o_n(1) \\
& = \theta\int_{\Omega}(\mu_1|z_n|^{4}+\mu_2|w_n|^{4})dx + 2\beta\int_{\Omega}|z_n|^{2}|w_n|^{2}dx\\
& \leq \theta\int_{\Omega}(\mu_1|z_n|^{4}+\mu_2|w_n|^{4})dx \\
& \leq \theta\max\{\mu_1,\mu_2\}\int_{\Omega}(|z_n|^{4}+|w_n|^{4})dx  \\
& \leq \frac{\theta\max\{\mu_1,\mu_2\}}{\mathcal{S}^2}\left(\int_{\Omega}(|\nabla z_n|^2+|\nabla w_n|^2)dx\right)^2.
\end{align*}
Let $\int_{\Omega}|\nabla z_n|^2dx \to k \geq 0, \int_{\Omega}|\nabla w_n|^2dx \to l \geq 0$ as $n \to \infty$. We have
\begin{align}
k + l \leq \frac{\theta\max\{\mu_1,\mu_2\}}{\mathcal{S}^2}(k+l)^2.
\end{align}
If $k + l > 0$, one gets $k+l \geq \frac{\mathcal{S}^2}{\theta\max\{\mu_1,\mu_2\}}$.
This shows that
$$
E_\theta(z_n,w_n) = \frac{1}{4}(k + l) + o_n(1) \geq \frac{1}{4\theta\max\{\mu_1,\mu_2\}}\mathcal{S}^2 + o_n(1).
$$
Then, using the Br\'{e}zis-Lieb Lemma (see \cite{BL}) one gets
$$
E_\theta(u_n,v_n) = E_\theta(u_\theta,v_\theta) + E_\theta(z_n,w_n) + o_n(1) \geq E_\theta(u_\theta,v_\theta) + \frac{1}{4\theta\max\{\mu_1,\mu_2\}}\mathcal{S}^2 + o_n(1),
$$
which contradicts to
$$
E_\theta(u_n,v_n) \to m_\theta < E_\theta(u_\theta,v_\theta) + \frac{1}{4\theta\max\{\mu_1,\mu_2\}}\mathcal{S}^2.
$$
Thus $k = l = 0$.

The fact that $k = l = 0$ implies $z_n \to 0, w_n \to 0$ strongly in $H_0^1(\Omega)$. This in turn yields to the fact that
$$
E_\theta(u_\theta,v_\theta) = \lim_{n \to \infty}E_\theta(u_n,v_n) = m_\theta.
$$
Hence $(u_\theta,v_\theta) \notin \mathcal{P}_\theta$. This proves that the claim holds true. If $(u_\theta,v_\theta) \in \partial\mathcal{P}_\theta$, then $E_\theta(u_\theta,v_\theta) = m_\theta$. Otherwise, $(u_\theta,v_\theta) \notin \overline{\mathcal{P}_\theta}$ and we aim to show that $E_\theta(u_\theta,v_\theta) \leq m_\theta$. Using \eqref{eq of zn}, \eqref{eq of wn} and \eqref{eq of zn 2}, \eqref{eq of wn 2} when $\beta > 0$ and $\beta < 0$ respectively, we obtain $E_\theta(z_n,w_n) \geq o_n(1)$. Then from
$$
m_\theta \leftarrow E_\theta(u_n,v_n) = E_\theta(u_\theta,v_\theta) + E_\theta(z_n,w_n) + o_n(1) \geq E_\theta(u_\theta,v_\theta) + o_n(1),
$$
one derives that $E_\theta(u_\theta,v_\theta) \leq m_\theta$.

By Lagrange multiplier principle, $(u_\theta,v_\theta)$ satisfies
\begin{equation}
\left\{
\begin{aligned}
-\Delta u_\theta = \lambda_{u_\theta} u_\theta^+ + \theta\mu_1 (u_\theta^+)^3 + \theta\beta u_\theta^+ (v_\theta^+)^2, \\
-\Delta v_\theta = \lambda_{v_\theta} v_\theta^+ + \theta\mu_2 (v_\theta^+)^3 + \theta\beta (u_\theta^+)^2 v_\theta^+,
\end{aligned}
\right. \nonumber
\end{equation}
when $\beta > 0$ and
\begin{equation}
\left\{
\begin{aligned}
-\Delta u_\theta = \lambda_{u_\theta} u_\theta^+ + \theta\mu_1 (u_\theta^+)^3 + \beta u_\theta^+ (v_\theta^+)^2, \\
-\Delta v_\theta = \lambda_{v_\theta} v_\theta^+ + \theta\mu_2 (v_\theta^+)^3 + \beta (u_\theta^+)^2 v_\theta^+,
\end{aligned}
\right. \nonumber
\end{equation}
when $\beta < 0$. Now we prove that $u_\theta,v_\theta$ are positive. We address the case $\beta > 0$ and the proof is similar when $\beta < 0$. Multiplying $u_\theta^-$ and integrating on $\Omega$ for the first equation yield that
$$
\int_\Omega|\nabla u_\theta^-|^2dx = 0,
$$
implying that $u_\theta \geq 0$. By strong maximum principle, $u_\theta > 0$. Similarly, we have $v_\theta > 0$. Thus $(u_\theta,v_\theta)$ solves \eqref{eq theta}.

Finally, we show that \eqref{inequality} and \eqref{inequality 2} hold true. In fact, by the Pohozaev identity, we have
\begin{align}
\int_{\Omega}(|\nabla u_\theta|^2+|\nabla v_\theta|^2)dx & - \frac{1}{2}\int_{\partial\Omega}(|\nabla u_\theta|^2+|\nabla v_\theta|^2)\sigma\cdot nd\sigma \nonumber \\
& = \frac{\theta N}{4}\int_{\Omega}(\mu_1 |u_\theta|^4 + 2\beta |u_\theta|^2 |v_\theta|^2 + \mu_2 |v_\theta|^4)dx,
\end{align}
when $\beta > 0$ and
\begin{align}
\int_{\Omega}(|\nabla u_\theta|^2+|\nabla v_\theta|^2)dx & - \frac{1}{2}\int_{\partial\Omega}(|\nabla u_\theta|^2+|\nabla v_\theta|^2)\sigma\cdot nd\sigma \nonumber \\
& = \frac{\theta N}{4}\int_{\Omega}(\mu_1 |u_\theta|^4 + \mu_2 |v_\theta|^4)dx + \frac{\beta N}{2}\int_{\Omega}|u_\theta|^2 |v_\theta|^2dx,
\end{align}
when $\beta < 0$. Since $\Omega$ is star-shaped with respect to $0$, $\sigma\cdot n > 0$. Hence we obtain \eqref{inequality} when $\beta > 0$ and \eqref{inequality 2} when $\beta < 0$.
\qed\vskip 5pt

\begin{lemma} \label{continuity}
Under the assumptions of Lemma \ref{M-P geometry}, we have $\lim_{\theta \to 1^-}m_\theta = m_1$.
\end{lemma}

\noindent\textbf{Proof.  } For any $(u,v) \in S_{\vec{c}}^+$, we have $E_\theta(u,v) \geq E(u,v)$ when $\theta < 1$. This shows that $\displaystyle\liminf_{\theta \to 1^-}m_\theta \geq m_1$. It is sufficient to prove that $\displaystyle\limsup_{\theta \to 1^-}m_\theta \leq m_1$. By the definition of $m_1$, for any $\epsilon > 0$, we can take a $\gamma_0 \in \Gamma$ such that
$$
\sup_{t \in [0,1]}E(\gamma_0(t)) < m_1 + \epsilon.
$$
For any $\theta_n \to 1^-$, we have
\begin{align}
m_{\theta_n} = & \inf_{\gamma \in \Gamma}\sup_{t \in [0,1]}E_{\theta_n}(\gamma(t)) \nonumber \\
\leq & \sup_{t \in [0,1]}E_{\theta_n}(\gamma_0(t)) \nonumber \\
= & \sup_{t \in [0,1]}E(\gamma_0(t)) + o_n(1) \nonumber \\
< & m_1 + \epsilon + o_n(1).
\end{align}
By the arbitrariness of $\epsilon$ one gets
$$
\lim_{n \to \infty}m_{\theta_n} \leq m_1,
$$
implying that
$$
\limsup_{\theta \to 1^-}m_\theta \leq m_1.
$$
This completes the proof.
\qed\vskip 5pt

\begin{corollary}[Existence of a bounded (PS) sequence] \label{bounded (PS) sequence}
Under the assumptions of Lemma \ref{M-P geometry}, there exists a bounded sequence $\{(u_n,v_n)\} \subset S_{\vec{c}}^+$ such that
$$
\lim_{n \to \infty}E(u_n,v_n) \leq m_1 \text{ and that } E'|_{S_{\vec{c}}^+}(u_n,v_n) \to 0 \text{ as } n \to \infty.
$$
Furthermore, if $N = 3$, then $\displaystyle\lim_{n \to \infty}E(u_n,v_n) = m_1$.
\end{corollary}

\noindent\textbf{Proof.} By Proposition \ref{Positive solutions for almost every theta}, we can take $\theta_n \to 1^-$ and $(u_n,v_n) = (u_{\theta_n},v_{\theta_n})$ solving \eqref{eq theta} when $\beta > 0$ and solving \eqref{eq theta 2} when $\beta < 0$. We firstly aim to show that $u_n, v_n$ are bounded in $H_0^1(\Omega)$. Note that $m_\theta$ is non-increasing on $[1-\epsilon,1]$. Hence, $m_\theta \leq m_1$, implying that when $\beta > 0$,
$$
\frac{1}{2}\int_{\Omega}(|\nabla u_n|^{2}+|\nabla v_n|^{2})dx - \frac{\theta_n}{4}\int_{\Omega}(\mu_1 u_n^4 + 2\beta u_n^2 v_n^2 + \mu_2 v_n^4)dx \leq m_1.
$$
Then using \eqref{inequality} we have
$$
\left(\frac{1}{2}-\frac{1}{N}\right)\int_{\Omega}(|\nabla u_n|^{2}+|\nabla v_n|^{2})dx \leq m_1,
$$
implying that $u_n, v_n$ are bounded in $H_0^1(\Omega)$. Quite similarly, when $\beta < 0$, one gets the boundedness of $u_n, v_n$ in $H_0^1(\Omega)$.

Next we prove that
$$
\lim_{n \to \infty}E(u_n,v_n) \leq m_1 \text{ and that } E'|_{S_{\vec{c}}^+}(u_n,v_n) \to 0 \text{ as } n \to \infty.
$$
We treat the case that $\beta > 0$. Following similar arguments, we can complete the proof when $\beta < 0$. By Proposition \ref{Positive solutions for almost every theta}, $E_{\theta_n}(u_n,v_n) \leq m_{\theta_n}$. Up to a subsequence, by Lemma \ref{continuity} we have
\begin{align}
\lim_{n \to \infty}E(u_n,v_n) & = \lim_{n \to \infty}\left(E_{\theta_n}(u_n,v_n) + \frac{\theta_n-1}{4}\int_\Omega (\mu_1 u_n^4 + 2\beta u_n^2 v_n^2 + \mu_2 v_n^4)dx\right) \nonumber \\
& = \lim_{n \to \infty}E_{\theta_n}(u_n,v_n) \leq \lim_{n \to \infty}m_{\theta_n} = m_1.
\end{align}
Moreover, as $n \to \infty$,
$$
\partial_uE|_{S_{\vec{c}}^+}(u_n,v_n) = \partial_uE_{\theta_n}|_{S_{\vec{c}}^+}(u_n,v_n) + (\theta_n-1)(\mu_1 u_n^3 + \beta u_n v_n^2) \to 0,
$$
$$
\partial_vE|_{S_{\vec{c}}^+}(u_n,v_n) = \partial_vE_{\theta_n}|_{S_{\vec{c}}^+}(u_n,v_n) + (\theta_n-1)(\mu_2 v_n^3 + \beta u_n^2 v_n) \to 0.
$$

Finally, if $N = 3$, by Proposition \ref{Positive solutions for almost every theta}, $E(u_n,v_n) = m_{\theta_n}$. Then by Lemma \ref{continuity}, one gets
$$
\lim_{n \to \infty}E(u_n,v_n) = \lim_{n \to \infty}E_{\theta_n}(u_n,v_n) = \lim_{n \to \infty}m_{\theta_n} = m_1.
$$
The proof is complete.
\qed\vskip 5pt

\section{Estimation of the M-P level when $N = 4$} \label{estimate}

Let $\xi \in C_0^\infty(\Omega)$ be the radial function, such that $\xi(x) \equiv 1$ for $0 \leq |x| \leq R$, $0 \leq \xi(x) \leq 1$ for $R \leq
|x| \leq 2R$, $\xi(x) \equiv 0$ for $|x| \geq 2R$, where $B_{2R} \subset \Omega$ and $|\nabla \xi| \leq 2/R$. Take $v_\epsilon = \xi U_\epsilon$ where
$$
U_\epsilon = \frac{2\sqrt{2}\epsilon}{\epsilon^2+|x|^2}.
$$
\begin{lemma}[] \label{lem es1}
If $N = 4$, then we have, as $\epsilon \to 0^+$,
\begin{equation}\label{es-1}
			\begin{aligned}
				\int_{\Omega} |\nabla v_\varepsilon|^2dx =\mathcal{S}^{2}+O(\varepsilon^{2}), \quad  	\int_{\Omega} |v_\varepsilon|^{4}dx =\mathcal{S}^{2}+O(\varepsilon^4).
			\end{aligned}
		\end{equation}
		\begin{equation}\label{es-2}
			\int_{\Omega}|v_\varepsilon|^2dx = d\varepsilon^2|\ln \varepsilon| +O(\varepsilon^2),
		\end{equation}
where $d$ is a positive constant, and as $\epsilon \to 0^+$,
		\begin{equation}\label{es-6}
			\int_{\Omega}|v_\varepsilon|^pdx \sim \epsilon^{4-p}, 2 < p < 4.
		\end{equation}
\end{lemma}

\begin{lemma}\label{lem es4}
If $N = 4$, then we have, for any positive $\varphi \in C(\Omega)\cap H^1(\Omega)$ as $\epsilon \to 0^+$,
\begin{equation}\label{es-7}
				\int_{\Omega}\varphi v_\varepsilon dx \leq 4\sqrt{2}\sup_{B_\delta}\varphi\omega_4R^2\epsilon + o(\epsilon),
		\end{equation}
		\begin{equation}\label{es-8}
			\int_{\Omega}\varphi v_\varepsilon^{3} dx \geq 4\sqrt{2}\inf_{B_\delta}\varphi\omega_4\epsilon + o(\epsilon),
		\end{equation}
where $\omega_4$ denotes the area of the unit sphere surface.
\end{lemma}

\textit{Proof.  } Direct computations yield that
\begin{align}
\int_{\Omega}\varphi v_\varepsilon dx =& 2\sqrt{2}\int_{B_{2R}}\varphi \xi\frac{\epsilon}{\epsilon^2+|x|^2}dx \nonumber \\
\leq& 2\sqrt{2}\sup_{B_\delta}\varphi \epsilon^{3} \int_{B_{\frac{2R}{\epsilon}}}\frac{1}{1+|x|^2}dx \nonumber \\
=& 2\sqrt{2}\sup_{B_\delta}\varphi \omega_4\epsilon^{3} \int_{0}^{\frac{2R}{\epsilon}}\frac{ r^{3}}{1+r^2}dr  \nonumber \\
\leq& 4\sqrt{2}\sup_{B_\delta}\varphi \omega_4R^2\epsilon + o(\epsilon),
\end{align}
\begin{align}
\int_{\Omega}\varphi v_\varepsilon^{3} dx =& 16\sqrt{2}\int_{B_{2R}}\varphi \xi^{3} \left(\frac{\epsilon}{\epsilon^2+|x|^2}\right)^{3}dx \nonumber \\
\geq& 16\sqrt{2}\inf_{B_\delta}\varphi \epsilon \int_{B_{\frac{R}{\epsilon}}}\left(\frac{1}{1+|x|^2}\right)^{3}dx \nonumber \\
=& 16\sqrt{2}\inf_{B_\delta}\varphi \omega_4\epsilon \int_{0}^{\frac{R}{\epsilon}}\left(\frac{r}{1+r^2}\right)^{3}dr  \nonumber \\
\geq& 4\sqrt{2}\inf_{B_\delta}\varphi \omega_4\epsilon + o(\epsilon).
\end{align}
\qed\vskip 5pt

\begin{theorem}[Estimation of the M-P level]\label{Estimate of the M-P level}
Assume that $N = 4$ and that \eqref{eq pre for m-p of E} holds for some $\rho > \sqrt{(c_1+c_2)\lambda_1}$. Let $(u,v) \in \mathcal{P}_1 \cap (C(\Omega) \times C(\Omega))$ be a positive solution of \eqref{eq1.1}, where $\mathcal{P}_1$ is defined by Lemma \ref{lem equal of m}. Then
\begin{align}
& m_1 < E(u,v)+ \frac{1}{4\max\{\mu_1,\mu_2\}}\mathcal{S}^2 \ \text{if} \ \beta < 0, \nonumber \\
& m_1 < E(u,v) + \frac{1}{2}\frac{1}{\sqrt{(\mu_1-\mu_2)^2+4\beta^2}+\mu_1+\mu_2}\mathcal{S}^2 \ \text{if} \ \beta > 0.
\end{align}
\end{theorem}

\noindent\textbf{Proof.}  Let $z_{\epsilon,s} = u + sv_\epsilon, w_{\epsilon,t} = v + tv_\epsilon, s, t \geq 0$. Then $z_{\epsilon,s}, w_{\epsilon,t} > 0$ in $\Omega$. If $w \in H_0^1(\Omega)$, then $w$ can be viewed as a function in $H^1(\mathbb{R}^N)$ by defining $w(x) = 0$ for all $x \notin \Omega$. Define
$$Z_{\epsilon,s} = \tau z_{\epsilon,s}(\tau x), \quad W_{\epsilon,t} = \eta w_{\epsilon,t}(\eta x).$$
By taking $$\tau = \frac{\|z_{\epsilon,s}\|_{L^2(\Omega)}}{\sqrt{c_1}} \geq 1, \quad \eta = \frac{\|w_{\epsilon,t}\|_{L^2(\Omega)}}{\sqrt{c_2}} \geq 1,$$
we obtain $Z_{\epsilon,s}, W_{\epsilon,t} \in H_0^1(\Omega)$ and $\|Z_{\epsilon,s}\|_{L^2(\Omega)}^2 = c_1$, $\|W_{\epsilon,t}\|_{L^2(\Omega)}^2 = c_2$. Let
$$\overline{Z}_{k} = k^{2}Z_{\epsilon,\widehat{s}}(kx), \overline{W}_{k} = k^{2}W_{\epsilon,\widehat{t}}(kx), k \geq 1,$$
where $\widehat{s}, \widehat{t}$ will be determined later. Set
\begin{align}
&\psi(k) \nonumber \\
&= E(\overline{Z}_{k},\overline{W}_{k}) \nonumber \\
& =\frac{k^2}{2}\int_{\Omega}(|\nabla Z_{\epsilon,\widehat{s}}|^2 + |\nabla W_{\epsilon,\widehat{t}}|^2dx) \nonumber \\
&\quad - \frac{k^{4}}{4}\int_{\Omega}(\mu_1|Z_{\epsilon,\widehat{s}}|^4 + 2\beta|Z_{\epsilon,\widehat{s}}|^2|W_{\epsilon,\widehat{t}}|^2 + \mu_2|W_{\epsilon,\widehat{t}}|^4)dx, k \geq 1.
\end{align}
Then,
\begin{align}
\psi'(k) = & k\int_{\Omega}(|\nabla Z_{\epsilon,\widehat{s}}|^2 + |\nabla W_{\epsilon,\widehat{t}}|^2dx) \nonumber \\
& - k^{3}\int_{\Omega}(\mu_1|Z_{\epsilon,\widehat{s}}|^4 + 2\beta|Z_{\epsilon,\widehat{s}}|^2|W_{\epsilon,\widehat{t}}|^2 + \mu_2|W_{\epsilon,\widehat{t}}|^4)dx.
\end{align}
By Lemma \ref{lem es1},
$$
\int_{\Omega}|\nabla Z_{\epsilon,\widehat{s}}|^2dx = \int_{\Omega}|\nabla u|^2dx + \widehat{s}^2\mathcal{S}^{2} + o_\epsilon(1),
$$
$$
\int_{\Omega}|Z_{\epsilon,\widehat{s}}|^{4}dx = \int_{\Omega}|u|^{4}dx + \widehat{s}^{4}\mathcal{S}^{2} + o_\epsilon(1),
$$
$$
\int_{\Omega}|Z_{\epsilon,\widehat{s}}|^{2}dx = \int_{\Omega}|u|^{2}dx + o_\epsilon(1),
$$
$$
\int_{\Omega}|\nabla W_{\epsilon,\widehat{t}}|^2dx = \int_{\Omega}|\nabla v|^2dx + \widehat{t}^2\mathcal{S}^{2} + o_\epsilon(1),
$$
$$
\int_{\Omega}|W_{\epsilon,\widehat{t}}|^{4}dx = \int_{\Omega}|v|^{4}dx + \widehat{t}^{4}\mathcal{S}^{2} + o_\epsilon(1),
$$
$$
\int_{\Omega}|W_{\epsilon,\widehat{t}}|^{2}dx = \int_{\Omega}|v|^{2}dx + o_\epsilon(1),
$$
one can choose $\max\{\widehat{s},\widehat{t}\}$ large such that $\psi'(k) < 0$ for all $k > 1$. Hence, $E(\overline{Z}_{k},\overline{W}_{k}) \leq E(Z_{\epsilon,\widehat{s}},W_{\epsilon,\widehat{t}})$. Furthermore, it not difficult to see that $\psi(k) \to -\infty$ as $k \to \infty$.

Note that direct computation yields
\begin{align} \label{1}
&E(Z_{\epsilon,s},W_{\epsilon,t}) \nonumber \\
& = E(z_{\epsilon,s},w_{\epsilon,t}) \nonumber \\
& \quad + \frac{\beta}{2}\int_{\Omega}(z_{\epsilon,s}(x))^2(w_{\epsilon,t}(x))^2dx - \frac{\beta \tau^2\eta^2}{2} \int_{\Omega}(z_{\epsilon,s}(\mu x))^2(w_{\epsilon,t}(\eta x))^2dx,
\end{align}
\begin{align} \label{2}
E(z_{\epsilon,s},w_{\epsilon,t}) = E(u,v) & + s\int_{\Omega}(\nabla u \nabla v_\epsilon - \mu_1u^3v_\epsilon - \beta uv^2v_\epsilon)dx  \nonumber \\
& + t\int_{\Omega}(\nabla v \nabla v_\epsilon - \mu_1v^3v_\epsilon - \beta u^2vv_\epsilon)dx \nonumber \\
& + \frac{1}{2}(s^2+t^2)\int_{\Omega}|\nabla v_\epsilon|^2dx - \frac{1}{4}(\mu_1s^4+2\beta s^2t^2 + \mu_2t^4)\int_{\Omega}|v_\epsilon|^4dx \nonumber \\
& - \frac{1}{2}\int_{\Omega}(3\mu_1s^2u^2 + \beta t^2u^2 + 3\mu_2t^2v^2 + \beta s^2v^2 + 4\beta stuv)v_\epsilon^2dx \nonumber \\
& - \int_{\Omega}(\mu_1s^3u + \mu_2t^3v +2t^2su + 2ts^2v)v_\epsilon^3dx \nonumber \\
= E(u,v) & + s\lambda_u\int_{\Omega}uv_\epsilon dx + t\lambda_v\int_{\Omega}vv_\epsilon dx \nonumber \\
& + \frac{1}{4}\left(2(s^2+t^2)-(\mu_1s^4+2\beta s^2t^2 + \mu_2t^4)\right)\mathcal{S}^2 \nonumber \\
& - \int_{\Omega}(\mu_1s^3u + \mu_2t^3v +2t^2su + 2ts^2v)v_\epsilon^3dx + o(\epsilon),
\end{align}
where in the computation of \eqref{2}, we used Lemma \ref{lem es1} and the fact that $(u,v)$ is a positive solution of \eqref{eq1.1}. By Lemma \ref{lem es4}, there exists $K_1 = K_1(\widehat{s},\widehat{t},\lambda_u,\lambda_v,u,v) > 0$ and $K_2 = K_2(\mu_1,\mu_2,\widehat{s},\widehat{t},u,v) > 0$ such that
\begin{align} \label{3}
& \left|s\lambda_u\int_{\Omega}uv_\epsilon dx + t\lambda_v\int_{\Omega}vv_\epsilon dx\right| \leq K_1R^2\epsilon + o(\epsilon), \nonumber \\
& - \int_{\Omega}(\mu_1s^3u + \mu_2t^3v +2t^2su + 2ts^2v)v_\epsilon^3dx \leq - K_2\epsilon + o(\epsilon).
\end{align}
Furthermore, we estimate
\begin{align}
\int_{\Omega}(z_{\epsilon,s}(x))^2(w_{\epsilon,t}(x))^2dx - \tau^2\eta^2 \int_{\Omega}(z_{\epsilon,s}(\tau x))^2(w_{\epsilon,t}(\eta x))^2dx.
\end{align}
Note that
\begin{align}
& \int_{\Omega}(z_{\epsilon,s}(x))^2(w_{\epsilon,t}(x))^2dx - \tau^2\eta^2 \int_{\Omega}(z_{\epsilon,s}(\tau x))^2(w_{\epsilon,t}(\eta x))^2dx \nonumber \\
= & (1-\tau^2\eta^2)\int_{\Omega}(z_{\epsilon,s}(x))^2(w_{\epsilon,t}(x))^2dx \nonumber \\
& + \tau^2\eta^2\int_{\Omega}\left((z_{\epsilon,s}(x))^2(w_{\epsilon,t}(x))^2 -(z_{\epsilon,s}(\tau x))^2(w_{\epsilon,t}(x))^2\right)dx \nonumber \\
& + \tau^2\eta^2\int_{\Omega}\left((z_{\epsilon,s}(\tau x))^2(w_{\epsilon,t}(x))^2 - (z_{\epsilon,s}(\tau x))^2(w_{\epsilon,t}(\eta x))^2\right)dx \nonumber \\
= & (1-\tau^2\eta^2)\int_{\Omega}(z_{\epsilon,s}(x))^2(w_{\epsilon,t}(x))^2dx \nonumber \\
& + 2\tau^2\eta^2(1-\tau)\int_{\Omega}z_{\epsilon,s}(x)x\cdot \nabla z_{\epsilon,s}(x) (w_{\epsilon,t}(x))^2dx + o(1-\tau) \nonumber \\
& + 2\tau^2\eta^2(1-\eta)\int_{\Omega}(z_{\epsilon,s}(\tau x))^2w_{\epsilon,t}(x)x\cdot \nabla w_{\epsilon,t}(x)dx + o(1-\eta).
\end{align}
For small $\epsilon$ and bounded $s, t$, we can assume $\tau, \eta \in [1,2]$. Since
\begin{align}
\int_{\Omega}|U_\epsilon|^4dx = \int_{\Omega}|U_1|^4dx, \quad \int_{\Omega}|x \cdot \nabla U_\epsilon|^4dx = \int_{\Omega}|x \cdot \nabla U_1|^4dx,
\end{align}
it is easy to verify that
$$\left|\int_{\Omega}z_{\epsilon,s}(x)x\cdot \nabla z_{\epsilon,s}(x) (w_{\epsilon,t}(x))^2dx\right| \leq \frac{C}{R},$$
$$
\left|\int_{\Omega}(z_{\epsilon,s}(\tau x))^2w_{\epsilon,t}(x)x\cdot \nabla w_{\epsilon,t}(x)dx\right| \leq \frac{C}{R},
$$
and that
$$\left|\int_{\Omega}(z_{\epsilon,s}(\tau x))^2(w_{\epsilon,t}(\eta x))^2dx\right| \leq C,$$ where $C$ is independent of $\epsilon$ and $R$. From the definitions of $\eta, \mu$ and Lemma \ref{lem es4}, one gets
$$
1-\mu = \frac{1-\mu^2}{1+\mu} \lesssim R^2\epsilon + o(\epsilon),
$$
$$
1-\eta = \frac{1-\eta^2}{1+\eta} \lesssim R^2\epsilon + o(\epsilon),
$$
$$
1-\mu^2\eta^2 \lesssim R^2\epsilon + o(\epsilon).
$$
Hence, there exists $K_3 = K_3(\widehat{s},\widehat{t},c_1,c_2,u,v,\beta) > 0$ such that
\begin{align} \label{4}
& \left|\frac{\beta}{2}\int_{\Omega}(z_{\epsilon,s}(x))^2(w_{\epsilon,t}(x))^2dx - \frac{\beta\tau^2\eta^2}{2} \int_{\Omega}(z_{\epsilon,s}(\tau x))^2(w_{\epsilon,t}(\eta x))^2dx\right| \nonumber \\
\leq & K_3(R+R^2)\epsilon + o(\epsilon).
\end{align}
By \eqref{1}, \eqref{2}, \eqref{3}, \eqref{4}, we have
\begin{align}
E(Z_{\epsilon,s},W_{\epsilon,t}) \leq E(u,v) & + \frac{1}{4}\left(2(s^2+t^2)-(\mu_1s^4+2\beta s^2t^2 + \mu_2t^4)\right)\mathcal{S}^2 \nonumber \\
& + (K_1R^2+K_3(R+R^2)-K_2)\epsilon + o(\epsilon).
\end{align}
By taking $R$ small such that $K_1R^2+K_3(R+R^2)-K_2 < 0$ and $\epsilon$ small enough, one gets
\begin{align}
E(Z_{\epsilon,s},W_{\epsilon,t}) < E(u,v) + \frac{1}{4}\left(2(s^2+t^2) -(\mu_1s^4+2\beta s^2t^2 + \mu_2t^4)\right)\mathcal{S}^2 - \delta_\epsilon ,
\end{align}
for some $\delta_\epsilon > 0$. Next, we consider three cases.

\vskip0.1in
\noindent \textbf{Case 1:} When $\beta < 0$ and $\mu_2 \leq \mu_1$, we choose $t \equiv 0$. Note that $2s^2 - \mu_1s^4 \leq 1/\mu_1$. Hence,
\begin{align}
\sup_{s > 0}E(Z_{\epsilon,s},W_{\epsilon,0}) < E(u,v) + \frac{1}{4\mu_1}\mathcal{S}^2.
\end{align}
Define $\widetilde{\gamma}(r):= (Z_{\epsilon,2r\widehat{s}},W_{\epsilon,0})$ for $r \in [0,1/2]$ and $$\widetilde{\gamma}(r):= (\overline{Z}_{2(r-1/2)k_0+1},\overline{W}_{2(r-1/2)k_0+1})$$ for $r \in [1/2,1]$ where $k_0$ is large such that $E(\overline{Z}_{k_0+1},\overline{W}_{k_0+1}) < E(u,v)$. Then, we have $\widetilde{\gamma} \in \tilde{\Gamma}$ and that $\displaystyle\sup_{r \in [0,1]}E(\widetilde{\gamma}) < E(u,v) + \frac{1}{4\mu_1}\mathcal{S}^2$, implying that $\tilde{m}_1 < E(u,v)+ \frac{1}{4\mu_1}\mathcal{S}^2$.

\vskip0.1in

\noindent\textbf{Case 2:} When $\beta < 0$ and $\mu_2 > \mu_1$, we choose $s \equiv 0$. Note that $2t^2 - \mu_2t^4 \leq 1/\mu_2$. Hence,
\begin{align}
\sup_{t > 0}E(Z_{\epsilon,0},W_{\epsilon,t}) < E(u,v) + \frac{1}{4\mu_2}\mathcal{S}^2.
\end{align}
Quite similar to the Case 1, we have $\tilde{m}_1 < E(u,v) + \frac{1}{4\mu_2}\mathcal{S}^2$.

\vskip0.1in

\noindent\textbf{Case 3:} When $\beta > 0$,  we  choose $\epsilon_0 s^2 = t^2$ where
$$
\epsilon_0 = \frac{\sqrt{(\mu_1-\mu_2)^2+4\beta^2}-\mu_1+\mu_2}{2\beta}.
$$
such that $2s^2t^2 = \epsilon_0 s^4 + 1/\epsilon_0 t^4$ and $\mu_1 + \epsilon_0\beta = \mu_2 + \beta/\epsilon_0$. Note that
\begin{align}
& 2(s^2+t^2) -(\mu_1s^4+2\beta s^2t^2 + \mu_2t^4) \nonumber \\
= & 2(1+\epsilon_0)s^2 - (\mu_1+\epsilon_0\beta)(1+\epsilon_0^2)s^4 \nonumber \\
\leq & \frac{1}{\sqrt{(\mu_1-\mu_2)^2+4\beta^2}+\mu_1+\mu_2}\frac{(1+\epsilon_0)^2}{1+\epsilon_0^2} \nonumber \\
\leq & \frac{2}{\sqrt{(\mu_1-\mu_2)^2+4\beta^2}+\mu_1+\mu_2}.
\end{align}
Hence,
\begin{align}
\sup_{s > 0,t = \sqrt{\epsilon_0}s}E(Z_{\epsilon,s},W_{\epsilon,t}) < E(u,v) + \frac12\frac{1}{\sqrt{(\mu_1-\mu_2)^2+4\beta^2}+\mu_1+\mu_2}\mathcal{S}^2.
\end{align}
Define $\widetilde{\gamma}(r):= (Z_{\epsilon,2r\widehat{s}},W_{\epsilon,2r\widehat{t}})$ for $r \in [0,1/2]$ and $$\widetilde{\gamma}(r):= (\overline{Z}_{2(r-1/2)k_0+1},\overline{W}_{2(r-1/2)k_0+1})$$ for $r \in [1/2,1]$ where $k_0$ is large such that $E(\overline{Z}_{k_0+1},\overline{W}_{k_0+1}) < E(u,v)$. Then $\widetilde{\gamma} \in \tilde{\Gamma}$ and $\displaystyle\sup_{r \in [0,1]}E(\widetilde{\gamma}) < E(u,v) + \frac{1}{2}\frac{1}{\sqrt{(\mu_1-\mu_2)^2+4\beta^2}+\mu_1+\mu_2}\mathcal{S}^2$, implying that $$\tilde{m}_1 < E(u,v) + \frac12\frac{1}{\sqrt{(\mu_1-\mu_2)^2+4\beta^2}+\mu_1+\mu_2}\mathcal{S}^2.$$
Then by Lemma \ref{lem equal of m}, we complete the proof.
\qed\vskip 5pt

\section{Proof of Theorem \ref{thmB.5}} \label{second solution}

\textbf{Proof of Theorem \ref{thmB.5}.  } Note that all the assumptions in Lemma \ref{M-P geometry} hold. Thus there is a M-P structure. By Corollary \ref{bounded (PS) sequence}, there is a bounded (PS) sequence $\{(u_n,v_n)\}$ at the level less than or equal to $m_1$.

If $N = 3$, which is a Sobolev subcritical case, using Sobolev compact embedding, by a standard process, we can prove that $u_n \to u, v_n \to v$ strongly in $H_0^1(\Omega)$. Furthermore, by Corollary \ref{bounded (PS) sequence}, $E(u_n,v_n) \to m_1$ as $n \to \infty$, implying that $E(u,v) = m_1$ and $(u,v)$ is of M-P type. By taking $(\tilde{u}_{\vec{c}},\tilde{v}_{\vec{c}}) = (u,v)$ we complete the proof when $N = 3$.

Next we focus on the Sobolev critical case ($N = 4$). Up to a subsequence, we assume that
\begin{equation}
\begin{aligned}
&u_n \rightharpoonup u, \quad v_n \rightharpoonup v, \ \text{ weakly in } H_0^1(\Omega), \\
&u_n \rightharpoonup u, \quad v_n \rightharpoonup v, \ \text{ weakly in } L^{4}(\Omega),\\
&u_n \to u, \quad v_n \to v, \ \text{ strongly in } L^p(\Omega) \text{ for } 2 \leq p<4,\\
&u_n \to u, \quad v_n \to v, \ \text{ almost everywhere in } \Omega.
\end{aligned}
\end{equation}
It can be verified that $(u,v) \in S_{\vec{c}}^+$ is a solution of \eqref{eq1.1}. Furthermore, $u,v$ are positive. If $(u,v) \neq (u_{\vec{c}},v_{\vec{c}})$, we complete the proof by taking $(\tilde{u}_{\vec{c}},\tilde{v}_{\vec{c}}) = (u,v)$. If not, assume that $(u,v) \equiv (u_{\vec{c}},v_{\vec{c}})$ and we aim to find a contradiction.

\

\noindent\textbf{Step 1:} There exists $\widehat{\epsilon}>0$ such that $(u_{\vec{c}},v_{\vec{c}})$ is in $\mathcal{P}_\theta \cap \mathcal{B}_{\rho}$ for all $\theta \in [1-\widehat{\epsilon},1)$ where $\mathcal{P}_\theta$ is given by Lemma \ref{lem equal of m} and $\rho$ is given by Lemma \ref{M-P geometry}, or the proof is complete.

Setting
$$
\nu_\theta = \inf_{(z,w) \in \mathcal{P}_\theta \cap \mathcal{B}_{\rho}}E(z,w) > -\infty,
$$
we will show that $\nu_\theta$ can be achieved by a $(z_0,w_0) \in \mathcal{P}_\theta \cap \mathcal{B}_{\rho}$ to complete the proof of this step. We give the proof when $\beta > 0$ and the case when $\beta < 0$ can be proved after similar discussions. Note that when $\beta > 0$,
\begin{align}
& \inf_{(z,w) \in \partial(\mathcal{P}_\theta \cap \mathcal{B}_{\rho})}E(z,w) \nonumber \\
=& \inf_{(z,w) \in \partial(\mathcal{P}_\theta \cap \mathcal{B}_{\rho})}\left(E_\theta(z,w) + \frac{\theta-1}{4}\int_\Omega(\mu_1|z^+|^4+2\beta|z^+|^2|w^+|^2+\mu_2|w^+|^4))dx\right) \nonumber \\
=& \inf_{(z,w) \in (\partial\mathcal{P}_\theta \cap \overline{\mathcal{B}_{\rho}}) \cup (\overline{\mathcal{P}_\theta} \cap \partial \mathcal{B}_{\rho})}E_\theta(z,w) + O(\theta-1) \nonumber \\
>& E_\theta(\sqrt{c_1}e_1,\sqrt{c_2}e_1) + \delta + O(\theta-1) \nonumber \\
\geq & E(\sqrt{c_1}e_1,\sqrt{c_2}e_1) + \delta + O(\theta-1),
\end{align}
for all $\theta \in [1-\epsilon,1)$ where $\epsilon$ and $\delta$ are given by Lemma \ref{M-P geometry}. There exists $\epsilon' > 0$ such that $|O(\theta-1)| < \delta/2$ for all $\theta \in [1-\epsilon',1)$. Choosing $\widehat{\epsilon} = \min\{\epsilon,\epsilon'\}$ we obtain
$$
\inf_{(z,w) \in \partial(\mathcal{P}_\theta \cap \mathcal{B}_{\rho})}E(z,w) > E(\sqrt{c_1}e_1,\sqrt{c_2}e_1) +\frac{\delta}{2}
$$
for all $\theta \in [1-\widehat{\epsilon},1)$. Hence, we can take a minimizing sequence $(z_n,w_n)$ of $\nu_\theta, \theta  \in [1-\widehat{\epsilon},1)$ staying in
$$
\Pi_\tau := \{(z,w) \in \mathcal{P}_\theta \cap \mathcal{B}_{\rho}: \text{dist}((z,w), \partial(\mathcal{P}_\theta \cap \mathcal{B}_{\rho})) \geq \tau\}
$$
for some $\tau > 0$ small enough with $\displaystyle\sup_n E_\theta(z_n,w_n) < m_\theta$. Obviously, $z_n, w_n$ are bounded in $H_0^1(\Omega)$. Up to a subsequence, we assume that
\begin{equation}
\begin{aligned}
&z_n \rightharpoonup z_0, \quad w_n \rightharpoonup w_0, \ \text{ weakly in } H_0^1(\Omega), \\
&z_n \rightharpoonup z_0, \quad w_n \rightharpoonup w_0, \ \text{ weakly in } L^{4}(\Omega),\\
&z_n \to z_0, \quad w_n \to w_0, \ \text{ strongly in } L^p(\Omega) \text{ for } 2\leq p<4,\\
&z_n \to z_0, \quad w_n \to w_0, \ \text{ almost everywhere in } \Omega.
\end{aligned}
\end{equation}
If $\int_\Omega(|\nabla(z_n-z_0)|^2+|\nabla(w_n-w_0)|^2)dx \to 0$ as $n \to \infty$, we complete the proof of this step. Hence, we assume that
\begin{align} \label{eq ass}
\liminf_{n \to \infty}\int_\Omega(|\nabla(z_n-z_0)|^2+|\nabla(w_n-w_0)|^2)dx > 0.
\end{align}
By the lower semicontinuity of the norm, one gets $(z_0,w_0) \in \mathcal{B}_{\rho}$. Similar to the discussions in Proposition \ref{M-P geometry}, we can show that $E(z_0,w_0) \leq \nu_\theta$. Now we aim to show that $(z_0,w_0) \in \mathcal{P}_\theta$ and thus $E(z_0,w_0) \geq \nu_\theta$. Let $\gamma(t) = (z_0,w_0) + t(z_n-z_0,w_n-w_0)$. By the Br\'{e}zis-Lieb Lemma (see \cite{BL}) again one can get
\begin{align}
& \int_\Omega(|\nabla(z_n-z_0)|^2+|\nabla(w_n-w_0)|^2)dx = \nonumber \\
& \int_\Omega(\mu_1|(z_n-z_0)^+|^4+2\beta|(z_n-z_0)^+|^2|(w_n-w_0)^+|^2+\mu_2|(w_n-w_0)^+|^4)dx\nonumber \\
&\quad  + o_n(1).
\end{align}
By \eqref{eq ass} we have
$$
\frac{\int_\Omega(|\nabla(z_n-z_0)|^2+|\nabla(w_n-w_0)|^2)dx}{\theta\int_\Omega(\mu_1|(z_n-z_0)^+|^4+2\beta|(z_n-z_0)^+|^2|(w_n-w_0)^+|^2+\mu_2|(w_n-w_0)^+|^4))dx}$$
 $$\to \frac{1}{\theta} > 1.
$$
Hence for large $n$,
\begin{align}
& \sup_{t \in [0,1]}\bigg(\frac{t^2}{2}\int_\Omega(|\nabla(z_n-z_0)|^2+|\nabla(w_n-w_0)|^2)dx - \nonumber \\
& \frac{\theta t^4}{4}\int_\Omega(\mu_1|(z_n-z_0)^+|^4+2\beta|(z_n-z_0)^+|^2|(w_n-w_0)^+|^2+\mu_2|(w_n-w_0)^+|^4))dx\bigg) \nonumber \\
= & \frac{1}{2}\int_\Omega(|\nabla(z_n-z_0)|^2+|\nabla(w_n-w_0)|^2)dx - \nonumber \\ &\frac{\theta}{4}\int_\Omega(\mu_1|(z_n-z_0)^+|^4+2\beta|(z_n-z_0)^+|^2|(w_n-w_0)^+|^2+\mu_2|(w_n-w_0)^+|^4)dx.
\end{align}
Then for $n$ large enough, using Br\'{e}zis-Lieb Lemma we have
\begin{align}
& E_\theta(\gamma(t)) \nonumber \\
=& o_n(1) + E_\theta(z_0,w_0) + \frac{t^2}{2}\int_\Omega(|\nabla(z_n-z_0)|^2+|\nabla(w_n-w_0)|^2)dx - \nonumber \\
& \frac{\theta t^4}{4}\int_\Omega(\mu_1|(z_n-z_0)^+|^4+2\beta|(z_n-z_0)^+|^2|(w_n-w_0)^+|^2+\mu_2|(w_n-w_0)^+|^4)dx  \nonumber \\
\leq& o_n(1) + E_\theta(z_0,w_0) + \frac{1}{2}\int_\Omega(|\nabla(z_n-z_0)|^2+|\nabla(w_n-w_0)|^2)dx - \nonumber \\
& \frac{\theta}{4}\int_\Omega(\mu_1|(z_n-z_0)^+|^4+2\beta|(z_n-z_0)^+|^2|(w_n-w_0)^+|^2+\mu_2|(w_n-w_0)^+|^4)dx \nonumber \\
\leq& E_\theta(z_0,w_0) + E_\theta(z_n-z_0,w_n-w_0) + o_n(1) \nonumber \\
=& E_\theta(z_n,w_n) + o_n(1), \nonumber
\end{align}
where $o_n(1)$ is independent of $t \in [0,1]$. Thus for large $n$, $$\sup_{t \in [0,1]}E_\theta(\gamma(t)) \leq \sup_nE_\theta(z_n,w_n) < m_\theta.$$ Hence, $(z_0,w_0) \in \mathcal{P}_\theta$. Together with $(z_0,w_0) \in \mathcal{B}_{\rho}$, we have $E(z_0,w_0) \geq \nu_\theta$. Then, one gets $E(z_0,w_0) = \nu_\theta$.

Note that $(z_0,w_0)$ is a positive solution of \eqref{eq1.1}-\eqref{eq1.2}. If $(u_{\vec{c}},v_{\vec{c}}) \neq (z_0,w_0)$, we find two positive solutions and the proof is complete by taking $(\tilde{u}_{\vec{c}},\tilde{v}_{\vec{c}}) = (z_0,w_0)$. From now on, we always assume that there exists $\widehat{\epsilon}>0$ such that $(u_{\vec{c}},v_{\vec{c}})$ is in $\mathcal{P}_\theta \cap \mathcal{B}_{\rho}$ for all $\theta \in [1-\widehat{\epsilon},1)$.

\

\noindent\textbf{Step 2:} Completion.  Let $z_n = u_n - u_{\vec{c}}$, $w_n = v_n - v_{\vec{c}}$. Using the Br\'{e}zis-Lieb Lemma (see \cite{BL}) one gets
$$
\int_{\Omega}|\nabla z_n|^2dx = \int_{\Omega}(\mu_1|z_n|^{4}+\beta|z_n|^{2}|w_n|^{2})dx + o_n(1),
$$
$$
\int_{\Omega}|\nabla w_n|^2dx = \int_{\Omega}(\mu_2|w_n|^{4}+\beta|z_n|^{2}|w_n|^{2})dx + o_n(1),
$$
$$
E(u_n,v_n) = E(u_{\vec{c}},v_{\vec{c}}) + E(z_n,w_n) + o_n(1).
$$
Assume that $\int_{\Omega}|\nabla z_n|^2dx \to k \geq 0, \int_{\Omega}|\nabla w_n|^2dx \to l\geq 0$, as $n \to \infty$.

\vskip0.1in
On the one hand, the case of $k + l = 0$ is impossible. Suppose on the contrary that $k = l = 0$, then $z_n \to 0$, $w_n \to 0$ strongly in $H_0^1(\Omega)$. Hence $u_n \to u_{\vec{c}}$, $v_n \to v_{\vec{c}}$ strongly in $H_0^1(\Omega)$. By Proposition \ref{Positive solutions for almost every theta}, $E(u_n,v_n) = m_{\theta_n}+o_n(1)$ or $(u_n,v_n) \notin \overline{\mathcal{P}_{\theta_n}}$. If the former holds, one gets
$$
E(u_{\vec{c}},v_{\vec{c}}) = \lim_{n \to \infty}m_{\theta_n} \geq E(\sqrt{c_1}e_1,\sqrt{c_2}e_1)+\delta \geq E(u_{\vec{c}},v_{\vec{c}})+\delta,
$$
which is impossible. The latter is in a contradiction with $(u_{\vec{c}},v_{\vec{c}}) \in \mathcal{P}_{\theta_n}$ for large $n$ that is shown in Step 1. In fact, set $(u^t_n,v^t_n) = (u_{\vec{c}} + t(u_n-u_{\vec{c}}),v_{\vec{c}} + t(v_n-v_{\vec{c}})), t \in [0,1]$. Then $(u^t_n,v^t_n) \to (u_{\vec{c}},v_{\vec{c}})$ as $n \to \infty$ for all $t \in [0,1]$. However, there exists $t_0^n \in (0,1)$ such that $(u_n^{t_0^n},v_n^{t_0^n}) \in \partial \mathcal{P}_{\theta_n}$, which implies that
$$
E(u_n^{t_0^n},v_n^{t_0^n}) = E_{\theta_n}(u_n^{t_0^n},v_n^{t_0^n}) + o_n(1) = m_{\theta_n} + o_n(1) \geq E(u_{\vec{c}},v_{\vec{c}})+\delta + o_n(1).
$$
This is a contradiction.

On the other hand, the case of $k + l > 0$ is impossible. If not, similar to the proof of Proposition \ref{Positive solutions for almost every theta}, we can find a contradiction with Theorem \ref{Estimate of the M-P level}. The proof is complete.
\qed\vskip 5pt

\section{Proof of Theorem \ref{thmps}} \label{secps}

This section is devoted to the proof of Theorem \ref{thmps}. Recall that $\lambda_2 := \lambda_2(\Omega)$ is the second eigenvalue of $-\Delta$ on $\Omega$ with Dirichlet boundary condition and $e_2$ is the corresponding unit eigenfunction, $e_2^+ = \max\{e_2,0\}$ and $e_2^- = e_2 - e_2^+$.

\

\noindent\textbf{Proof of Theorem \ref{thmps}}
Using \eqref{eq beta < 0}, \eqref{eq < 0 2} and the fact that
$$
E\left(\sqrt{c_1}\frac{e_2^+}{\|e_2^+\|_{L^2(\Omega)}}, -\sqrt{c_2}\frac{e_2^-}{\|e_2^-\|_{L^2(\Omega)}}\right) < \frac{\lambda_2(c_1+c_2)}{2} \ \text{when} \ \beta < 0,
$$
we have
$$
E\left(\sqrt{c_1}\frac{e_2^+}{\|e_2^+\|_{L^2(\Omega)}}, -\sqrt{c_2}\frac{e_2^-}{\|e_2^-\|_{L^2(\Omega)}}\right) < \inf_{(u,v) \in \partial\mathcal{B}_{\rho}\cap S_{\vec{c}}^+}E(u,v).
$$
Then following the steps of proving Theorems \ref{thmB.2} and \ref{thmB.5}, one can show that \eqref{eq1.1}-\eqref{eq1.2} has a positive local minimizer $(u_{1,\beta},v_{1,\beta})$ and the second positive solution $(u_{2,\beta},v_{2,\beta})$.

The asymptotic behavior of $(u_{1,\beta},v_{1,\beta})$ as $\beta \to -\infty$ can be proved following similar discussions to \cite[Section 5]{NTV2} and we omit the details here.

Next we study the asymptotic behavior of $(u_{2,\beta},v_{2,\beta})$. Similar to Lemma \ref{eq pre for m-p}, we can define $m_{1,\beta}$:
\begin{align}
m_{1,\beta} := \inf_{\gamma \in \Gamma}\sup_{t \in [0,1]}E(\gamma(t)) > & E\left(\sqrt{c_1}\frac{e_2^+}{\|e_2^+\|_{L^2(\Omega)}}, -\sqrt{c_2}\frac{e_2^-}{\|e_2^-\|_{L^2(\Omega)}}\right) \nonumber \\
= & \max\left\{E\left(\sqrt{c_1}\frac{e_2^+}{\|e_2^+\|_{L^2(\Omega)}}, -\sqrt{c_2}\frac{e_2^-}{\|e_2^-\|_{L^2(\Omega)}}\right),E_\theta(z,w)\right\},
\end{align}
where
$$
\Gamma := \left\{\gamma \in C([0,1],S_{\vec{c}}^+): \gamma(0) = \left(\sqrt{c_1}\frac{e_2^+}{\|e_2^+\|_{L^2(\Omega)}}, -\sqrt{c_2}\frac{e_2^-}{\|e_2^-\|_{L^2(\Omega)}}\right), \gamma(1) = (z,w)\right\},
$$
for some $(z,w) \in S_{\vec{c}}^+ \cap \mathcal{B}_{\rho}^c$. We may assume that $E(u_{2,\beta},v_{2,\beta}) \leq m_{1,\beta}$. Let $\mathcal{P}_{1,\beta} \subset S_{\vec{c}}^+$ be the maximal path connected branch containing $\left(\sqrt{c_1}\frac{e_2^+}{\|e_2^+\|_{L^2(\Omega)}}, -\sqrt{c_2}\frac{e_2^-}{\|e_2^-\|_{L^2(\Omega)}}\right)$ such that
$$
E(u,v) < m_{1,\beta}, \forall (u,v) \in \mathcal{P}_{1,\beta}.
$$
Further, by the proof of Theorem \ref{thmB.5}, we can assume that $(u_{1,\beta},v_{1,\beta}) \in \mathcal{P}_{1,\beta}$. Also, we know
$$
E(u_{1,\beta},v_{1,\beta}) \leq E\left(\sqrt{c_1}\frac{e_2^+}{\|e_2^+\|_{L^2(\Omega)}}, -\sqrt{c_2}\frac{e_2^-}{\|e_2^-\|_{L^2(\Omega)}}\right).
$$
Then, using a similar estimation of $m_{1,\beta}$ to Theorem \ref{Estimate of the M-P level}, one gets
\begin{align} \label{eq bound}
&E(u_{2,\beta},v_{2,\beta}) \leq m_{1,\beta}   < E(u_{1,\beta},v_{1,\beta}) + \frac{1}{4\max\{\mu_1,\mu_2\}}\mathcal{S}^2 \nonumber \\
& \leq E\left(\sqrt{c_1}\frac{e_2^+}{\|e_2^+\|_{L^2(\Omega)}}, -\sqrt{c_2}\frac{e_2^-}{\|e_2^-\|_{L^2(\Omega)}}\right) + \frac{1}{4\max\{\mu_1,\mu_2\}}\mathcal{S}^2 := M,
\end{align}
where $M$ is independent of $\beta < 0$. By the Pohozaev identity, we have
\begin{align}
	& ~~~~~~ \int_{\Omega}(|\nabla u_{2,\beta}|^2+|\nabla v_{2,\beta}|^2)dx \nonumber \\
	& \geq \int_{\Omega}(|\nabla u_{2,\beta}|^2+|\nabla v_{2,\beta}|^2)dx - \frac{1}{2}\int_{\partial\Omega}(|\nabla u_{2,\beta}|^2+|\nabla v_{2,\beta}|^2)\sigma\cdot nd\sigma \nonumber \\
	& = \frac{N}{4}\int_{\Omega}(\mu_1 |u_{2,\beta}|^4 + 2\beta |u_{2,\beta}|^2 |v_{2,\beta}|^2 + \mu_2 |v_{2,\beta}|^4)dx,
\end{align}
since $\sigma\cdot n > 0$ when $\Omega$ is star-shaped with respect to $0$. Then it holds
\begin{align} \label{eq m}
& M > E(u_{2,\beta},v_{2,\beta}) \nonumber \\
& = \frac12\int_{\Omega}(|\nabla u_{2,\beta}|^2+|\nabla v_{2,\beta}|^2)dx - \frac14\int_{\Omega}(\mu_1 |u_{2,\beta}|^4 + 2\beta |u_{2,\beta}|^2 |v_{2,\beta}|^2 + \mu_2 |v_{2,\beta}|^4)dx \nonumber \\
& \geq \left( \frac12-\frac1N\right) \int_{\Omega}(|\nabla u_{2,\beta}|^2+|\nabla v_{2,\beta}|^2)dx.
\end{align}
Hence, both $\{u_{2,\beta}\}_{\beta < 0}$ and $\{v_{2,\beta}\}_{\beta < 0}$ are uniformly bounded in $H_0^1(\Omega)$, and so in $L^4(\Omega)$. Moreover, using \eqref{eq m} one gets
\begin{align}
-2\beta\int_{\Omega}|u_{2,\beta}|^2 |v_{2,\beta}|^2dx < 4M - & 2\int_{\Omega}(|\nabla u_{2,\beta}|^2+|\nabla v_{2,\beta}|^2)dx \nonumber \\
+ & \int_{\Omega}(\mu_1 |u_{2,\beta}|^4 + \mu_2 |v_{2,\beta}|^4)dx < 2M_1,
\end{align}
for some $M_1$ which is independent of $\beta < 0$, thus $\{|\beta|\int_{\Omega}|u_{2,\beta}|^2 |v_{2,\beta}|^2dx\}_{\beta < 0}$ is uniformly bounded. Since $(u_{2,\beta},v_{2,\beta})$ solves \eqref{eq1.1}, by testing the first equation in \eqref{eq1.1} by $u_{2,\beta}$ and the second one by $v_{2,\beta}$, and using the previous estimates, we have,
\begin{align}
c_1|\eta_{2,\beta}| \leq & \left|\int_{\Omega}(|\nabla u_{2,\beta}|^2 - \mu_1|u_{2,\beta}|^4 - \beta|u_{2,\beta}|^2 |v_{2,\beta}|^2) dx\right| \nonumber \\
\leq & \int_{\Omega}(|\nabla u_{2,\beta}|^2 + \mu_1|u_{2,\beta}|^4 + |\beta||u_{2,\beta}|^2 |v_{2,\beta}|^2) dx \nonumber \\
\leq & \frac{2N}{N-2} M + \mu_1\mathcal{C}_N^4\left( \frac{2N}{N-2} M\right) ^2 + M_1,
\end{align}
\begin{align}
c_2|\xi_{2,\beta}| \leq & \left|\int_{\Omega}(|\nabla v_{2,\beta}|^2 - \mu_2|v_{2,\beta}|^4 - \beta|u_{2,\beta}|^2 |v_{2,\beta}|^2) dx\right| \nonumber \\
\leq & \int_{\Omega}(|\nabla v_{2,\beta}|^2 + \mu_2|v_{2,\beta}|^4 + |\beta||u_{2,\beta}|^2 |v_{2,\beta}|^2) dx \nonumber \\
\leq & \frac{2N}{N-2} M + \mu_2\mathcal{C}_N^4\left( \frac{2N}{N-2} M\right) ^2 + M_1.
\end{align}
Now using a Brezis-Kato-Moser type argument exactly as in \cite[pp. 1264-1265]{NTTV2} or \cite[The proof of Lemma 6.1]{Chen-Zou-2012}, we can obtain uniform $L^\infty$-bounds for $\{u_{2,\beta}\}_{\beta < 0}$ and $\{v_{2,\beta}\}_{\beta < 0}$. Using \cite[Theorems 1.3 and 1.5]{STTZ} (see also \cite{NTTV,SZ}), we know that these two sequences are uniformly bounded in $C^{0,\alpha}(\Omega)$ for every $0 < \alpha < 1$, and there exist $(u_\infty,v_\infty) \in C^{0,\alpha}(\Omega) \times C^{0,\alpha}(\Omega)$ with $u_\infty, v_\infty \geq 0$ in $\Omega$, and $\eta_2, \xi_2$ such that, up to subsequences, as $\beta \to -\infty$, it holds
\begin{align*}
u_{2,\beta} \to u_\infty, v_{2,\beta} \to v_\infty \text{ in } H_0^1(\Omega) \cap C^{0,\alpha}(\overline{\Omega}), \quad \eta_{2,\beta} \to \eta_2, \quad \xi_{2,\beta} \to \xi_2.
\end{align*}
So $u_{2,\beta} \to u_\infty, v_{2,\beta} \to v_\infty$ strongly in $L^2(\Omega)$ and then $\int_{\Omega}|u_\infty|^2dx = c_1$, $\int_{\Omega}|v_\infty|^2dx = c_2$. The fact that $|\beta|\int_{\Omega}|u_{2,\beta}|^2 |v_{2,\beta}|^2dx < M_1$ yields $\int_{\Omega}|u_{2,\beta}|^2 |v_{2,\beta}|^2dx \to 0$ as $\beta \to -\infty$. Hence, $u_\infty v_\infty = 0$ since $u_\infty, v_\infty$ are continuous. Then by \cite[Theorem 1.1]{DWZ}, we have
\begin{align*}
-\Delta (u_\infty - v_\infty) \geq \eta_2u_\infty - \xi_2v_\infty + \mu_1u_\infty^3 - \mu_2v_\infty^3  \text{ in } \Omega, \\
-\Delta (v_\infty - u_\infty) \geq \xi_2v_\infty - \eta_2u_\infty + \mu_2v_\infty^3 - \mu_1u_\infty^3 \text{ in } \Omega.
\end{align*}
By taking $Q_2 = u_\infty - v_\infty$ we complete the proof.
\qed\vskip 10pt




\begin{thebibliography}{99}
\bibitem{Akhmediev-Ankiewicz} N. Akhmediev, A. Ankiewicz, Partially coherent solitons on a finite background, Phys. Rev. Lett. 82 (1999) 2661-2664.
\bibitem{Ambrosetti-Colorado-2006} A. Ambrosetti, E. Colorado, Bound and ground states of coupled nonlinear Schr\"{o}dinger equations, C. R. Math. Acad. Sci. Paris 342 (2006) 453-458.
\bibitem{Ambrosetti-Colorado-2007} A. Ambrosetti, E. Colorado, Standing waves of some coupled nonlinear Schr\"{o}dinger equations, J. Lond. Math. Soc. 75 (2007) 67-82.
\bibitem{BCJS}  J. Borthwick, X. Chang, L. Jeanjean, N. Soave, Bounded Palais-Smale sequences with Morse type information for some constrained functionals, arXiv:2210.12626.
\bibitem{BL}  H. Br\'ezis, E.H. Lieb, A relation between pointwise convergence of functions and convergence of functionals, Proc. Amer. Math. Soc. 88 (3) (1983) 486-490.
\bibitem{BLZ} T. Bartsch, H.W. Li, W.M. Zou, Existence and asymptotic behavior of normalized ground states for Sobolev critical Schr\"{o}dinger systems, Calculus of Variations and Partial Differential Equations 62 (1) (2023) 1-34.
\bibitem{BJN}  T. Bartsch, L. Jeanjean, N. Soave, Normalized solutions for a system of coupled cubic Schr\"{o}dinger equations on $\mathbb{R}^3$, J. Math. Pures Appl. 106 (4) (2016) 583-614.
\bibitem{BN}  T. Bartsch, N. Soave, A natural constraint approach to normalized solutions of nonlinear Schr\"{o}dinger equations and systems, J. Funct. Anal. 272 (12) (2017), 4998-5037.
\bibitem{Bartsch-Wang}  T. Bartsch, Z.-Q. Wang, Note on ground states of nonlinear Schr\"{o}dinger systems, J. Partial Differ. Equ. 19 (2006) 200-207.
\bibitem{Bartsch-Wang-Wei}  T. Bartsch, Z.-Q. Wang, J. Wei, Bound states for a coupled Schr\"{o}dinger system, J. Fixed Point Theory Appl. 2 (2007) 353-367.
\bibitem{CHMS}  X.J. Chang, H. Hajaiej, Z.J. Ma, L.J. Song, Normalized ground states of the biharmonic nonlinear Schr\"{o}dinger equation with energy-critical exponent and combined nonlinearities, arXiv:2305.00327.
\bibitem{CL}  L.A. Caffarelli, F.-H. Lin, Singularly perturbed elliptic systems and multi-valued harmonic functions with free boundaries, J. Amer. Math. Soc. 21 (2008) 847-862.
\bibitem{CR}  L.A. Caffarelli, J.M. Roquejoffre, Uniform H\"{o}lder estimates in a class of elliptic systems and applications to singular limits in models for diffusion flames, Arch. Ration. Mech. Anal. 183 (2007) 457-487.
\bibitem{CTV}  M. Conti, S. Terracini, G. Verzini, Asymptotic estimates for the spatial segregation of competitive systems, Adv. Math. 195 (2005) 524-560.
\bibitem{Chen-Zou-2012}  Z.J. Chen, W.M. Zou, Positive least energy solutions and phase separation for coupled Schr\"{o}dinger equations with critical exponent, Arch. Ration. Mech. Anal. 205 (2012) 515-551.
\bibitem{DWZ}  E.N. Dancer, K.L. Wang, Z.T. Zhang, The limit equation for the Gross-Pitaevskii equations and S. Terracini's conjecture, J. Funct. Anal. 262 (3) (2012) 1087-1131.
\bibitem{EGBB}  B. Esry, C. Greene, J. Burke, J.  Bohn, Hartree-Fock theory for double condensates, Phys. Rev. Lett. 78 (1997) 3594-3597.
\bibitem{EHLS}  Amin Esfahani, H. Hajaiej, Y.M. Luo, L.J. Song, On the focusing fractional nonlinear Schr\"{o}dinger equation on the waveguide manifolds, arXiv:2305.19791.
\bibitem{HLS}  H. Hajaiej, Y.M. Luo, L.J. Song, On existence and stability results for normalized ground states of mass-subcritical biharmonic NLS on $\mathbb{R}^d \times \mathbb{T}^n$, arXiv:2212.00750.
\bibitem{HPS}  H. Hajaiej, E. Pacherie, L.J. Song, On the Number of Normalized Ground State Solutions for a class of Elliptic Equations with general nonlinearities and potentials, arXiv:2308.14599.
\bibitem{HS} H. Hajaiej, L.J. Song, A General and Unified Method to prove the Uniqueness of Ground State Solutions and the Existence/Non-existence, and Multiplicity of Normalized Solutions with applications to various NLS, arXiv:2208.11862.
\bibitem{HS2}  H. Hajaiej, L.J. Song, Strict Monotonicity of the global branch of solutions and Uniqueness of the corresponding normalized ground states for various classes of PDEs: Two general Methods with some examples, arXiv:2302.09681.
\bibitem{Lin-Wei}  T. Lin, J. Wei, Ground state of $N$ coupled nonlinear Schr\"{o}dinger equations in $\mathbb{R}^n, n \leq 3$, Commun. Math. Phys. 255 (2005) 629-653.
\bibitem{LZ}  H.W. Li, W.M. Zou, Normalized ground states for semilinear elliptic systems with critical and subcritical nonlinearities, J. Fixed Point Theory Appl. 23 (2021), Article number: 43.
\bibitem{Jean}  L. Jeanjean, On the existence of bounded Palais-Smale sequences and application to a Landesman-Lazer-type problem set on $\mathbb{R}^N$, Proc. Roy. Soc. Edinburgh Sect. A 129 (1999) 787-809.
\bibitem{JL}  L. Jeanjean, T.T. Le, Multiple normalized solutions for a Sobolev critical Schr\"{o}dinger equation, Math. Ann. 384 (2022) 101-134.
\bibitem{Maia-Pellacci-Squassina}  L. Maia, B. Pellacci, M. Squassina, Semiclassical states for weakly coupled nonlinear Schr\"{o}dinger systems, J. Eur. Math. Soc. 10 (2007) 47-71.
\bibitem{NTTV}  B. Noris, H. Tavares, S. Terracini, G. Verzini, Uniform H\"{o}lder bounds for nonlinear Schr\"{o}dinger systems with strong competition, Commun. Pure Appl. Math. 63 (2010) 267-302.
\bibitem{NTTV2}  B. Noris, H. Tavares, S. Terracini, G. Verzini, Convergence of minimax structures and continuation of critical points for singularly perturbed systems, J. Eur. Math. Soc. 14 (4) (2012) 1245-1273.
\bibitem{NTV2}  B. Noris, H. Tavares, G. Verzini, Normalized solutions for Nonlinear Schr\"{o}dinger systems on bounded domains, Nonlinearity 32 (2019).
\bibitem{Sirakov}  B. Sirakov, Least energy solitary waves for a system of nonlinear Schr\"{o}dinger equations in $\mathbb{R}^n$, Commun. Math. Phys. 271 (2007) 199-221.
\bibitem{Song} L.J. Song, Properties of the least action level, bifurcation phenomena and the existence of normalized solutions for a family of semi-linear elliptic equations without the hypothesis of autonomy, J. Differential Equations 315 (2022) 179-199.
\bibitem{Song2}  L.J. Song, Existence and orbital stability/instability of standing waves with prescribed mass for the $L^{2}$-supercritical NLS in bounded domains and exterior domains, Calc. Var. Partial Differential Equations 62 (6) (2023), Art. 176.
\bibitem{Song3}  L.J. Song, H. Hajaiej, Threshold for Existence, Non-existence and Multiplicity of positive solutions with prescribed mass for an NLS with a pure power nonlinearity in the exterior of a ball, arXiv:2209.06665.
\bibitem{Song4}  L.J. Song, H. Hajaiej, A New Method to prove the Existence, Non-existence, Multiplicity, Uniqueness, and Orbital Stability/Instability of standing waves for NLS with partial confinement, arXiv:2211.10058.
\bibitem{Song5}  L.J. Song, W.M. Zou, Two Positive Normalized Solutions on Star-shaped Bounded Domains to the Br\'ezis-Nirenberg Problem, I: Existence, in preparation.
\bibitem{STTZ}  N. Soave, H. Tavares, S. Terracini, A. Zilio, H\"{o}lder bounds and regularity of emerging free boundaries for strongly competing Schr\"{o}dinger equations with nontrivial grouping, Nonlinear Anal. 138 (2016) 388-427.
\bibitem{SZ}  N. Soave, A. Zilio, Uniform bounds for strongly competing systems: the optimal Lipschitz case, Arch. Ration. Mech. Anal. 218 (2) (2015) 647-697.
\bibitem{WW1}  J. Wei, T. Weth, Radial solutions and phase separation in a system of two coupled Schr\"{o}dinger equations, Arch. Ration. Mech. Anal. 190 (2008) 83-106.
\bibitem{WW2}  J. Wei, T. Weth, Asymptotic behaviour of solutions of planar elliptic systems with strong competition. Nonlinearity 21 (2008) 305-317.
\bibitem{WW} J. Wei, Y. Wu, Normalized solutions for Schr\"{o}dinger equations with critical Sobolev exponent and mixed nonlinearities, Journal of Functional Analysis 283 (6) (2022): 109574
\end{thebibliography}
\end{document}